\documentclass[a4paper]{article}
\usepackage{graphicx}
\usepackage{amssymb}
\usepackage{amsmath}
\usepackage{amsfonts}
\usepackage{varioref}
\usepackage[bookmarksnumbered]{hyperref}

\title{Convergence of penalty Robin--Robin domain decomposition methods for unilateral
 multibody contact problems of elasticity}

\author{Ivan~I.~Dyyak \footnote{Ivan Franko National University of Lviv,
Universytetska 1, Lviv, 79000, Ukraine, dyyak@franko.lviv.ua} \and
Ihor~I.~Prokopyshyn \footnote{Pidstryhach Institute for Applied
Problems of Mechanics and Mathematics, National Academy of
Sciences of Ukraine, Naukova 3-b, Lviv, 79060, Ukraine,
ihor84@gmail.com, Corresponding author } \and Ivan~A.~Prokopyshyn
\footnote {Ivan Franko National University of Lviv, Universytetska
1, Lviv, 79000, Ukraine, lviv.pi@gmail.com}}

\begin{document}
\maketitle

\begin{abstract}
The paper is devoted to the penalty Robin--Robin domain
decomposition methods (DDMs), proposed by us for the solution of
unilateral multibody contact problems of elasticity. These DDMs
are based on the penalty method for variational inequalities and
some stationary and nonstationary iterative methods for nonlinear
variational equations. The main result of the paper is that we
give the mathematical justification of proposed DDMs and prove
theorems on their convergence. We also investigate the numerical
efficiency of these methods using the finite element
approximations.


\end{abstract}

\section{Introduction}
\label{} The contact problems of elasticity are widely used in
many fields of science and engineering, especially in machine
science, structural mechanics, geology and biomechanics. The brief
overview of existing numerical and analytical methods for the
solution of contact problems can be found in \cite{PII_32,PII_33}.

Efficient approach for the solution of multibody contact problems
is the use of domain decomposition methods (DDMs).

DDMs are well developed for the solution of linear boundary value
problems, particularly for Poisson and linear elasticity problems
\cite{PII_31,PII_34,PII_35,PII_9,PII_28}. A Robin--Robin type
domain decomposition algorithm for linear Poisson boundary value
problems was introduced by P.~L.~Lions in \cite{PII_43}. Further,
Robin--Robin DDMs for linear elliptic boundary value problems were
investigated in works \cite{PII_44,PII_45}. An optimization based
domain decomposition methods for linear Poisson boundary value
problems were developed in \cite{PII_36,PII_37}.

The construction of DDMs for unilateral contact problems, which
are nonlinear, are much more complicated. Among the domain
decomposition methods for unilateral two-body contact problems
obtained on the continuous level, one should mention
Dirichlet--Neumann \cite {PII_2,PII_18,PII_38}, Neumann--Neumann
\cite {PII_3,PII_39} and optimization based \cite {PII_16}
iterative algorithms. A generalization of Lions' Robin--Robin
domain decomposition algorithm to a two-body contact problem was
proposed in works \cite {PII_40,PII_26,PII_41}. All of these
methods in each iteration require to solve a nonlinear one-sided
contact problem with a rigid body (Signorini problem) for one of
the bodies, and a linear elasticity problem with Neumann \cite
{PII_2,PII_18,PII_38} or Dirichlet \cite {PII_3,PII_39,PII_16}
boundary conditions on the possible contact area for the other
body, or require to solve nonlinear Signorini problems for both of
the bodies \cite {PII_40,PII_26,PII_41}. Moreover, to increase the
convergence rate of Neumann--Neumann and Robin--Robin algorithms,
it is recommended to perform an additional iteration, in which the
linear elasticity problems with Neumann boundary conditions have
to be solved for both of the bodies \cite {PII_3,PII_26}.

A domain decomposition method presented in work \cite {PII_17} for
two-body unilateral contact problem, is also obtained on
continuous level. It is based on the augmented Lagrangian
variational formulation and Uzawa block relaxation method. This
domain decomposition method in each iteration require to solve
linear elasticity problems with Robin boundary conditions for both
of the bodies.

On the contrary, DDMs can be constructed on the discrete level,
after a discretization of the corresponding continuous boundary
value problem. Among the discrete DDMs for unilateral contact
problems, one should mark out substructuring and FETI methods
\cite {PII_1,PII_5,PII_6,PII_42,PII_29}.

In works \cite {PII_23,PII_7,PII_8,PII_47} we proposed on the
continuous level a class of penalty parallel Robin--Robin type
domain decomposition methods for the solution of unilateral
multibody contact problems of elasticity. These methods are based
on the penalty method for variational inequalities and some
stationary and nonstationary iterative methods for nonlinear
variational equations. In each iteration of proposed DDMs we have
to solve in a parallel some linear variational equations in
subdomains, which correspond to linear elasticity problems with
Robin boundary conditions, prescribed on some subareas of the
possible contact zones. These DDMs do not require the solution of
nonlinear one-sided contact problems in each step.

The main result of this paper is that we prove theorems on the
convergence of proposed penalty Robin--Robin domain decomposition
methods. The paper is organized as follows. In section~2 the
classical formulation of the multibody contact problem in the form
of the system of second order elliptic partial differential
equations with inequality and equality constrains is given. In
section~3 we consider the variational formulations of this problem
in the form of convex minimization problem and in the form of
elliptic variational inequality at the closed convex set \cite
{PII_19,PII_20}. In section~4 we use the penalty method \cite
{PII_4,PII_21,PII_15} to reduce the variational inequality to an
unconstrained minimization problem, which is equivalent to a
nonlinear variational equation in the whole space. Later, we prove
a theorem on unique existence of a solution of the penalty
variational equation and a theorem on the strong convergence of
this solution to the solution of the original variational
inequality. In section~5 we consider stationary and nonstationary
iterative methods for the solution of abstract nonlinear
variational equations in reflexive Banach spaces. We prove
theorems on the convergence of these methods, and show that the
convergence rate of the stationary methods in some energy norm is
linear. We also formulate a theorem on stability of the stationary
iterative methods to the errors which may occur in each iteration.
In section~6 we present the parallel stationary and nonstationary
penalty Robin--Robin domain decomposition methods for the solution
of nonlinear penalty variational equations of unilateral multibody
contact problems. We prove a theorem on convergence of these
methods, and show that the convergence rate of the stationary
Robin--Robin methods in some energy norm is linear. In section~7
we perform the numerical analysis of proposed domain decomposition
methods using the finite element approximations. The penalty
parameter and the mesh refinement influence on the numerical
solution, as well as the dependence of the convergence rate of the
domain decomposition methods on the iterative parameters are
investigated. In conclusion section we summarize all results
presented in the paper.

\section{Formulation of unilateral multibody contact problem}

Introduce the Cartesian coordinate system $O{\kern 1pt} x_{1}
x_{2} x_{3}$ with basis vectors ${\bf e}_{1}, \, {\bf e}_{2}, \,
{\bf e}_{3}$, and consider the problem of frictionless unilateral
contact between $N$ elastic bodies $\Omega_{\alpha} \subset
{\mathbb R}^{3}$ with Lipschitz boundaries $\Gamma_{\alpha}
=\partial \Omega_{\alpha}$, $\alpha =1,2,...,N$ (Fig.~1). Denote
$\Omega =\bigcup_{\alpha=1}^{N}\Omega_{\alpha}$.

\begin{figure}[h]
 \center{
\includegraphics[bb=0mm 0mm 208mm 296mm, width=83.9mm, height=61.5mm,
viewport=3mm 4mm 205mm 292mm] {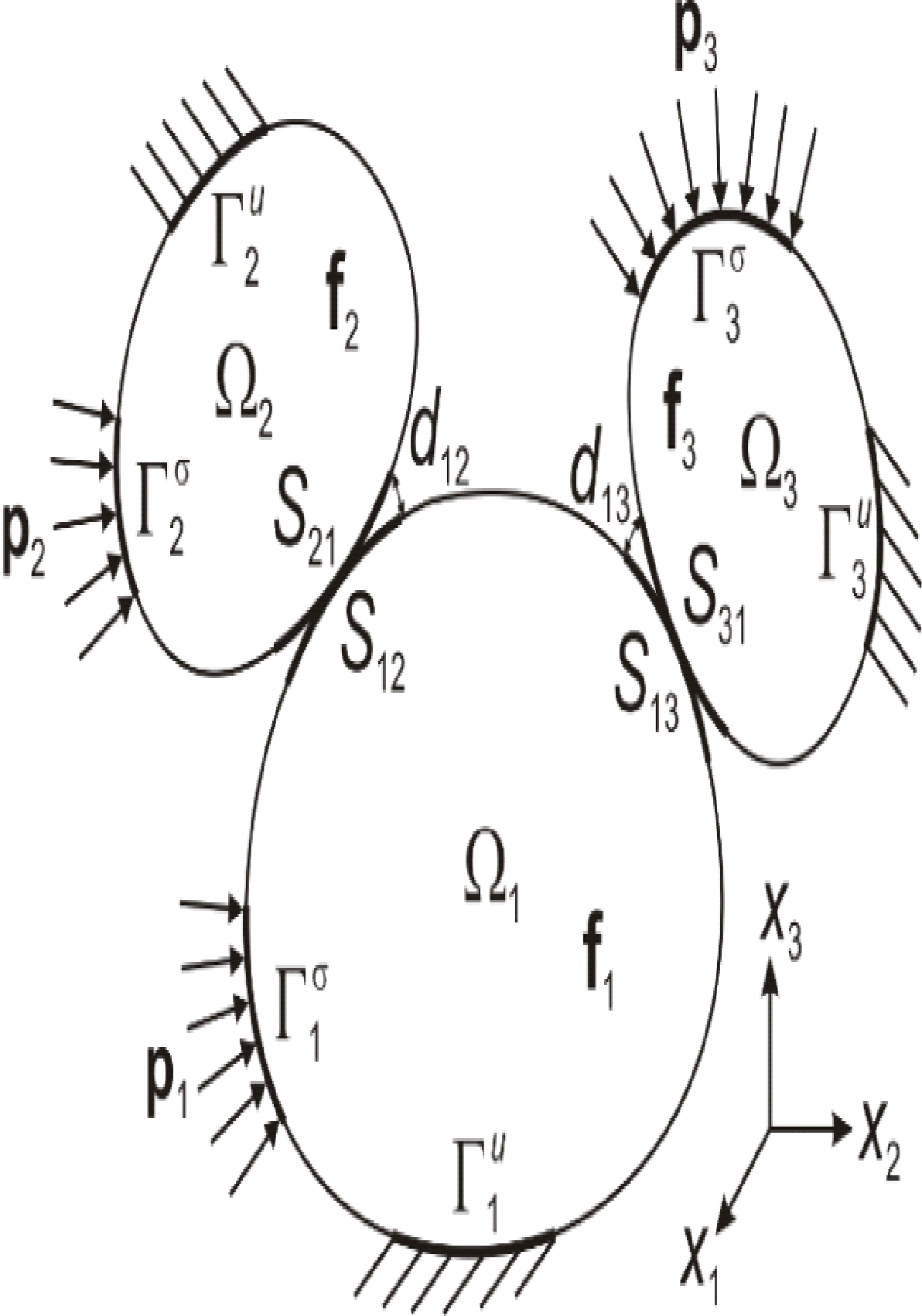}
 \\
 \textbf{Fig~1.} Unilateral contact between several elastic bodies 
 }
\end{figure}

The stress-strain state in point $ {\bf x}=(x_{1} ,x_{2} ,x_{3}
)^{{\rm T}}$ of each solid $\Omega_{\alpha}$ is descried by the
displacement vector ${\bf u}_{\, \alpha } ({\bf x})=u_{\alpha \,
i} ({\bf x})\, {\bf e}_{i} $, the symmetric tensor of strains
 $ {\hat{\pmb{\varepsilon}}}_{\alpha }=\varepsilon _{\alpha \, ij}\,
 {\bf e}_{i} \, {\bf e}_{j} $,
and the tensor of stresses $\hat{{\pmb \sigma }}_{\alpha } =\sigma
_{\alpha \, ij} \, {\bf e}_{i} \,
 {\bf e}_{j} $\,.
These quantities satisfy Cauchy relations, Hook's Law and the
equilibrium equations:
\begin{equation} \label{eq_PII_1}
\varepsilon_{\alpha \, ij}({\bf x})=\frac{1}{2}
\left(\frac{\partial u_{\alpha \, i}({\bf x})}{\partial x_{j}}
+\frac{\partial u_{\alpha \, j}({\bf x})}{\partial x_{i}} \right),
\,\,\, {\bf x} \in \Omega_{\alpha}, \,\,\, i,j=1,2,3,
\end{equation}
\begin{equation} \label{eq_PII_2}
\sigma_{\alpha \, ij}({\bf x})=\sum_{k,l=1}^{3}C_{\alpha \,ijkl}
({\bf x})\, \varepsilon_{\alpha \, kl} ({\bf x}), \,\,\, {\bf x}
\in \Omega_{\alpha}, \,\,\, i,j=1,2,3,
\end{equation}
\begin{equation} \label{eq_PII_3}
\sum_{j=1}^{3}\frac{\partial \sigma_{\alpha \, ij} ({\bf
x})}{\partial x_{j}} \, + f_{\alpha \, i} ({\bf x})=0, \, \,\,
{\bf x}\in \Omega_{\alpha} , \, \,\, i=1,2,3,
\end{equation}
where $f_{\alpha \,  i} ({\bf x})$ are the components of the
volume forces vector ${\bf f}_{\alpha } ({\bf x})=f_{\alpha \, i}
({\bf x})\, {\bf e}_{i} $.

The elastic coefficients $C_{\alpha \, ijkl} ({\bf x})$ are
measurable, symmetric, and uniformly elliptic with constants
$0<b\le d<\infty$:
\begin{equation} \label{eq_PII_4}
b\sum _{i,j=1}^{3}\varepsilon _{\alpha ij}^{2} ({\bf x})
 \le \sum _{i,j,k,l=1}^{3}C_{\alpha ijkl} ({\bf x})\, \varepsilon _{\alpha ij} ({\bf x})\,
\varepsilon _{\alpha kl} ({\bf x})\,  \le d\sum
_{k,l=1}^{3}\varepsilon _{\alpha kl}^{2} ({\bf x}) .
\end{equation}

Suppose that the boundary $\Gamma _{\alpha } $ of each solid
consists of three parts: $\Gamma _{\alpha }^{u} $, $\Gamma
_{\alpha }^{\sigma } $, $S_{\alpha } $, such that $\Gamma _{\alpha
} =\Gamma _{\alpha }^{u} \bigcup \Gamma _{\alpha }^{\sigma }
\bigcup S_{\alpha } $, $\Gamma _{\alpha }^{u} \bigcap \Gamma
_{\alpha }^{\sigma } \bigcap S_{\alpha } =\emptyset $, $\Gamma
_{\alpha }^{u} \ne \emptyset $, $\Gamma _{\alpha }^{u}
=\overline{\Gamma _{\alpha }^{u} }$, $S_{\alpha } \ne \emptyset $.
The boundary $S_{\alpha } =\bigcup _{\beta \in B_{\alpha }
}S_{\alpha \beta } $ is the possible contact area of the body
$\Omega _{\alpha } $ with the other bodies, $S_{\alpha \beta } $
is the possible contact area of the body $\Omega _{\alpha } $ with
the body $\Omega _{\beta } $, and $B_{\alpha } \subset
\left\{1,2,...,N\right\}$ is the set of the indices of all bodies
in contact with the body $\Omega_{\alpha}$.

On each boundary $\Gamma_{\alpha}$ let us introduce a local
orthonormal coordinate system ${\pmb \xi}_{\alpha},\, {\pmb \zeta
}_{\alpha},\, {\bf n}_{\,\alpha}$, where ${\bf n}_{\,\alpha}$ is
an outer unit normal to $\Gamma_{\alpha}$, and  ${\pmb
\xi}_{\alpha}$, ${\pmb \zeta}_{\alpha}$ are unit tangents. Then
the vectors of displacements and stresses on the boundary can be
written in the following way:
\[{\bf u}_{\alpha } ({\bf x})=u_{\alpha \, \xi } ({\bf x})\, {\pmb \xi }_{\alpha }
 +u_{\alpha \, \zeta } ({\bf x})\, {\pmb \zeta }_{\alpha } +u_{\alpha \, n}
({\bf x})\, {\bf n}_{\, \alpha }, \, \, {\bf x}\in \Gamma _{\alpha
} ,\]
\[{\pmb \sigma }_{\alpha } ({\bf x})={\hat{ \pmb{ \sigma }}}_{\alpha }
({\bf x})\cdot {\bf n}_{\, \alpha } =\sigma _{\alpha \, \xi }
({\bf x})\, {\pmb \xi }_{\, \alpha } +\sigma _{\alpha \, \zeta }
({\bf x})\, {\pmb \zeta }_{\alpha } +\sigma _{\alpha \, n} ({\bf
x})\,
 {\bf n}_{\, \alpha } , \, \, {\bf x}\in \Gamma _{\alpha } .\]

We assume that the surfaces $S_{\alpha \beta} \subset
\Gamma_{\alpha}$ and $S_{\beta \alpha} \subset \Gamma_{\beta}$ are
sufficiently close \cite {PII_19}. Therefore ${\bf n}_{\alpha}
({\bf x})\approx -{\bf n}_{\beta} ({\bf x'})$, where ${\bf
x'}=P({\bf x})\in S_{\beta \alpha}$ is an orthogonal projection of
point ${\bf x} \in S_{\alpha \beta}$ on the surface $S_{\beta
\alpha}$. We denote by $d_{\alpha \beta} ({\bf x})=\pm \left\|
{\bf x}-{\bf x'}\right\| _{2}=\pm
\sqrt{\sum_{j=1}^{3}\left(x_j-x'_j \right)^{2}}$ the distance in
${\mathbb R}^{3}$ between the bodies $\Omega_{\alpha}$ and
$\Omega_{\beta}$ before the deformation. The sign of $d_{\alpha
\beta} ({\bf x})$ depends on a statement of the specific problem.

On the part $\Gamma _{\alpha }^{u} $ the kinematical (Dirichlet)
boundary conditions are prescribed:
\begin{equation} \label{eq_PII_5}
{\bf u}_{\alpha } ({\bf x})={\bf z}_{\alpha } ({\bf x}), \, \,
{\bf x} \in \Gamma _{\alpha }^{u} ,
\end{equation}
and on the part $\Gamma _{\alpha }^{\sigma } $ we consider the
static (Neumann) boundary conditions
\begin{equation} \label{eq_PII_6}
{\pmb \sigma }_{\alpha } ({\bf x})={\bf p}_{\alpha } ({\bf x}), \,
\, {\bf x} \in \Gamma _{\alpha }^{\sigma } ,
\end{equation}
where ${\bf z}_{\alpha } ={z}_{\alpha \, \xi } ({\bf x})\, {\pmb
\xi }_{\alpha } +{z}_{\alpha \, \zeta } ({\bf x})\, {\pmb\zeta
}_{\alpha } +{z}_{\alpha \, n} ({\bf x})\, {\bf n}_{\alpha } $ and
$ {\bf p}_{\alpha } =p_{\alpha \, \xi } ({\bf x})\, {\pmb \xi
}_{\alpha } +p_{\alpha \, \zeta } ({\bf x})\, {\pmb \zeta
}_{\alpha } +p_{\alpha \, n} ({\bf x}) \, {\bf n}_{\alpha } $ are
given boundary displacements and stresses.

Further, for the simplicity of variational formulations and
proofs, we assume that all of the bodies are rigidly fixed on the
surface $\Gamma_{\alpha}^{u}$, i.e.
\begin{equation} \label{eq_PII_7}
{\bf z}_{\alpha} ({\bf x})=0, \,\, {\bf x} \in
\Gamma_{\alpha}^{u}.
\end{equation}
Note, that this assumption is not critical for the process of
numerical solution.

On the possible contact areas  $S_{\alpha \beta}$, $\alpha
=1,2,...,N$, $\beta \in B_{\alpha}$, the following unilateral
contact conditions hold:

\noindent absence of extension
\begin{equation} \label{eq_PII_8}
\sigma _{\alpha \, n} ({\bf x})=\sigma _{\beta \, n} ({\bf x'})\le
0,
\end{equation}
absence of friction
\begin{equation} \label{eq_PII_9}
\sigma _{\alpha \, \xi } ({\bf x})=\sigma _{\beta \, \xi } ({\bf
x'})=0, \,\, \sigma _{\alpha \, \zeta } ({\bf x})=\sigma _{\beta
\, \zeta } ({\bf x'})=0,
\end{equation}
mutual nonpenetration of the bodies
\begin{equation} \label{eq_PII_10}
u_{\alpha \, n} ({\bf x})+u_{\beta \, n} ({\bf x'})\le d_{\alpha
\beta } ({\bf x}),
\end{equation}
and contact alternative
\begin{equation} \label{eq_PII_11}
\left(u_{\alpha \, n} ({\bf x})+u_{\beta \, n} ({\bf
x'})-d_{\alpha \beta } ({\bf x})\, \right)\sigma _{\alpha \, n}
({\bf x})=0 \, ,
\end{equation}
where ${\bf x} \in S_{\alpha \beta } , \, \, {\bf x'}=P({\bf
x})\in S_{\beta \alpha }$.

The system of the second order partial differential equations
(\ref {eq_PII_1})~--~(\ref {eq_PII_3}) with the boundary
conditions (\ref {eq_PII_5})~--~(\ref {eq_PII_11}) is the
mathematical formulation of the frictionless unilateral multibody
contact problem of elasticity.

Note, that the contact problem (\ref {eq_PII_1})~--~(\ref
{eq_PII_3}), (\ref {eq_PII_5})~--~(\ref {eq_PII_11}) is nonlinear,
since the real contact areas are unknown.

\section{Variational formulation of the contact problem}

Let us consider the weak formulation of the contact problem (\ref
{eq_PII_1})~--~(\ref {eq_PII_3}), (\ref {eq_PII_5})~--~(\ref
{eq_PII_11}) in the form of variational inequality and convex
minimization problem. These variational formulations for the case
of the unilateral multibody contact problem were proposed in works
\cite {PII_19,PII_20}.

For each body $\Omega_{\alpha}$, $\alpha =1,2,...,N$, consider
Sobolev space $V_{\alpha } =[H^{1} (\Omega _{\alpha })]^{3}$ with
the scalar product $\left({\bf u}_{\alpha } ,{\bf v}_{\alpha }
\right)_{V_{\alpha } } =\sum _{i=1}^{3}\int _{\Omega _{\alpha }
}\left(u_{\alpha \, i} v_{\alpha \, i} +\sum
_{j=1}^{3}\frac{\partial u_{\alpha \, i} }{\partial x_{j} }
\frac{\partial v_{\alpha \, i} }{\partial x_{j} } \right)d\Omega
$,  \\ ${\bf u}_{\alpha} ,{\bf v}_{\alpha} \in V_{\alpha}$ and the
norm $\left\| {\bf u}_{\alpha } \right\| _{V_{\alpha } } =\sqrt{\,
\left({\bf u}_{\alpha } ,{\bf u}_{\alpha } \right)_{V_{\alpha } }
},\,\,{\bf u}_{\alpha } \in V_{\alpha }$.

Introduce the following closed subspace in $V_{\alpha}$:
\begin{equation} \label{eq_PII_12}
V_{\alpha }^{0} =\left\{\, {\bf u}_{\alpha }:\, \, \, \, \, {\bf
u}_{\alpha } \in V_{\alpha } \, ,\, \, \, \, {\rm{Tr}}_{\alpha
}^{u} ({\bf u}_{\alpha})=0\, \, \, \, {\rm on}\, \, \, \Gamma
_{\alpha }^{u} \, \right\},
\end{equation}
where ${\rm{Tr}}_{\alpha }^{u} :\, \, \, V_{\alpha } \to [H^{1/2}
(\Gamma _{\alpha }^{u} )]^{3}$ is surjective, linear and
continuous trace operator \cite {PII_22}. Space $V_{\alpha}^{0}$
is a Hilbert space with the same scalar product and norm as in
$V_{\alpha}$.

Consider the space $V_{0}$, which is the direct product of spaces
$V_{\alpha}^{0}$:
\begin{equation} \label{eq_PII_28}
V_{0}=V_{1}^{0} \times ...\times V_{N}^{0}=\left\{{\bf
u}=\left({\bf u}_{1},..., {\bf u}_{N} \right)^{{\rm T}}: \, \, \,
{\bf u}_{\alpha } \in V_{\alpha }^{0} ,\, \, \, \alpha
=1,2,...,N\right\},
\end{equation}
and define the scalar product and the norm in it: $\left({\bf
u},{\bf v}\right)_{V_{0} } =\sum _{\alpha =1}^{N}\left({\bf
u}_{\alpha } ,{\bf v}_{\alpha } \right)_{V_{\alpha } }  $,
$\left\| {\bf u}\right\| _{V_{0} } =\sqrt{\left({\bf u},{\bf
u}\right)_{V_{0} } } $, ${\bf u},{\bf v}\in V_{0} $. Note, that
the Hilbert space $V_{0} $ is a closed reflexive Banach space.

Now, let us introduce the closed convex set of all displacement
vectors in $V_{0}$ which satisfy the nonpenetration contact
conditions (\ref {eq_PII_10}):
\begin{equation} \label{eq_PII_13}
K=\left\{\, {\bf u}:\, \, \, \, {\bf u}\in V_{0} \, ,\, \, \, \,
u_{\alpha \, n} +u_{\beta \, n} \le d_{\alpha \beta } \, \, \,
{\rm on}\, \, \, S_{\alpha \beta } \, ,\, \, \, \, \left\{\alpha
,\, \beta \right\}\in Q\, \right\},
\end{equation}
where $Q=\left\{\, \left\{\alpha ,\beta \right\}:\, \, \, \alpha
\in \left\{1,2,...,N\right\},\, \, \, \beta \in B_{\alpha } \,
\right\}$ is the set of all possible unordered pairs of subscripts
of the bodies in contact with each other, and $d_{\alpha \beta }
\in H_{00}^{1/2} (\Xi_{\alpha } )$, $\left\{\alpha ,\, \beta
\right\}\in Q$, $\Xi_{\alpha }={\rm{int}}\, (\Gamma _{\alpha }
\backslash \Gamma _{\alpha }^{u} )$, $\alpha =1,2,...,N$.

The quantities $u_{\alpha \, n} $, $\alpha =1,2,...,N$, in (\ref
{eq_PII_13}) have to be understood in the following way
\[u_{\alpha \, n} ={\bf n}_{\alpha } \cdot {\rm{Tr}}_{\alpha }^{0} ({\bf u}_{\alpha}),\,\, {\bf u}_{\alpha } \in V_{\alpha }^{0} ,\]
where ${\rm{Tr}}_{\alpha }^{0}:\, \, \, V_{\alpha }^{0} \to
[H_{00}^{1/2} (\Xi _{\alpha } )]^{3} $ is surjective, linear and
continuous trace operator onto the surface $\Xi_{\alpha }
={\rm{int}}\, (\Gamma _{\alpha } \backslash \Gamma _{\alpha }^{u}
)$ \cite {PII_22}, and ${\bf n}_{\alpha } \in [L_{2} (\Xi _{\alpha
} )]^{3}$.

Note, that all equalities and inequalities in spaces $L_{2}$,
$H^{1/2}$, $H_{00}^{1/2}$ and $H^{1}$ hold almost everywhere.

Since the set $K$ is a closed convex subset of Hilbert space
$V_{0}$, it is weakly closed \cite {PII_4}.

In space $V_{0}$ consider a bilinear form $A\,(\bf{u},\bf{v})$,
such that $\frac{1}{2} A\,(\bf{u},\bf{u})$ represents the total
deformation energy of the system of bodies:
\begin{equation} \label{eq_PII_14}
A\, ({\bf u}, {\bf v})=\sum _{\alpha =1}^{N}a_{\alpha } ({\bf
u}_{\alpha } ,{\bf v}_{\alpha } ),\, \,\, {\bf u}, {\bf v}\in
V_{0} ,
\end{equation}

\begin{equation} \label{eq_PII_15}
a_{\alpha } ({\bf u}_{\alpha } ,{\bf v}_{\alpha } )=\int _{\Omega
_{\, \alpha } }{\hat{\pmb{\sigma}}}_{\alpha } ({\bf u}_{\alpha }
)\, :\,{\hat{\pmb{\varepsilon}}}_{\alpha } ({\bf v}_{\alpha } )\,
d\Omega  , \, \,\, {\bf u}_{\alpha } ,{\bf v}_{\alpha } \in
V_{\alpha }^{0} .
\end{equation}
Define in $V_{0}$ a linear form $L\,(\bf{v})$, which is equal to
the external forces work:
\begin{equation} \label{eq_PII_16}
L\, ({\bf v})=\sum _{\alpha =1}^{N} l_{\alpha } ({\bf v}_{\alpha }
) , \,\,\, {\bf v} \in V_{0} ,
\end{equation}

\begin{equation} \label{eq_PII_17}
l_{\alpha } ({\bf v}_{\alpha } )=\int _{\Omega _{\alpha } }{\bf
f}_{\alpha } \cdot {\bf v}_{\alpha } \, d\Omega  +\int _{\Gamma
_{\alpha }^{\, \sigma } }{\bf p}_{\alpha } \cdot {\rm{Tr}}_{\alpha
}^{0} (\, {\bf v}_{\alpha } ) \, dS, \,\,\, {\bf v}_{\alpha } \in
V_{\alpha }^{0} .
\end{equation}
where ${\bf f}_{\alpha} \in [L_{2}(\Omega_{\alpha})]^{3}$, ${\bf
p}_{\alpha} \in [L_{2}(\Gamma_{\alpha}^{\sigma})]^{3}$, $\alpha
=1,2,...,N$.

\textbf{Lemma~1.} \textit{If the boundaries} $\Gamma_{\alpha}
=\partial \Omega_{\alpha}$, $\alpha=1,2,...,N$, \textit{are
Lipschitz, $\Gamma _{\alpha}^{u} \ne \emptyset$},
$\Gamma_{\alpha}^{u} =\overline{\Gamma_{\alpha}^{u}}$, ${\bf
f}_{\alpha} \in [L_{2} (\Omega_{\alpha})]^{3}$, ${\bf p}_{\alpha}
\in [L_{2}(\Gamma_{\alpha}^{\sigma})]^{3}$, $C_{\alpha ijkl} \in
L_{\infty}(\Omega_{\alpha})$, $i,j,k,l=1,2,3$, $\alpha=1,2,...,N$,
\textit{and condition (\ref {eq_PII_4}) holds, then the bilinear
form $A$ is symmetric, continuous and coercive, and the linear
form $L$ is continuous, i.e.}

\begin{equation} \label{eq_PII_18}
\left(\forall {\bf u}, {\bf v} \in V_{0} \right)\left\{A\, ({\bf
u}, {\bf v})=A\, ({\bf v}, {\bf u})\right\},
\end{equation}

\begin{equation} \label{eq_PII_19}
\left(\exists M>0\right)\left(\forall {\bf u}, {\bf v}\in V_{0}
\right)\left\{\, \left|A\, ({\bf u}, {\bf v})\right|\le M \left\|
{\bf u}\right\| _{V_{0} } \left\| {\bf v}\right\| _{V_{0} }
\right\},
\end{equation}

\begin{equation} \label{eq_PII_20}
\left(\exists B>0\right)\left(\forall {\bf u} \in V_{0}
\right)\left\{A\, ({\bf u}, {\bf u})\ge B \left\| {\bf u}\right\|
_{V_{0}}^{2} \right\},
\end{equation}

\begin{equation} \label{eq_PII_21}
\left(\exists T>0\right)\left(\forall {\bf v}\in V_{0}
\right)\left\{\, \left|L\, (\bf{v})\right|\le T \left\| {\bf
v}\right\| _{V_{0}} \right\}.
\end{equation}

According to \cite {PII_19,PII_20}, the original contact problem
(\ref {eq_PII_1})~--~(\ref {eq_PII_3}), (\ref {eq_PII_5})~--~(\ref
{eq_PII_11}) has an alternative weak formulation as the convex
minimization problem of the quadratic functional on the set $K$:
\begin{equation} \label{eq_PII_22}
F({\bf u})=\frac{1}{2} A\, ({\bf u}, {\bf u})-L\, ({\bf u})\to
\mathop{\min }\limits_{{\bf u}\, \in K} .
\end{equation}
Using the general theory of variational inequalities \cite
{PII_4,PII_21,PII_22} the next theorem can be proved.

\textbf{Theorem~1.} \textit{Suppose that the conditions of
\textbf{\textit{Lemma~1}} hold and $d_{\alpha \beta} \in
H_{00}^{1/2}(\Xi_{\alpha})$. Then the minimization problem (\ref
{eq_PII_22}) has a unique solution on the convex set $K$, and this
problem is equivalent to the following variational inequality:}
\begin{equation} \label{eq_PII_23}
F'({\bf u},{\bf v}-{\bf u})=A\, ({\bf u},{\bf v}-{\bf u})-L\,
({\bf v}-{\bf u})\ge 0,\,\, \forall \, {\bf v}\in K .
\end{equation}

\section{Penalty variational formulation of the problem}

To obtain a minimization problem in the original space $V_{0}$, we
apply the penalty method \cite {PII_4,PII_21} to the convex
minimization problem (\ref {eq_PII_22}).

For the violation of nonpenetration conditions (\ref {eq_PII_10})
we use a penalty in the following form \cite {PII_15}:
\begin{equation} \label{eq_PII_24}
J_{\theta } ({\bf u})=\frac{1}{2\theta } \sum _{\left\{\alpha ,\,
\beta \right\}\, \in \, Q}\int _{S_{\alpha \beta }
}\left[\left(d_{\alpha \beta }
-u_{\alpha \, n}
-u_{\beta \, n}
\right)^{-}\right]^{2}\, dS,
\end{equation}
where $\theta >0$ is a penalty parameter, $y^{-} =\min \left\{0,\,
\, y\right\}$.

Let us consider the following minimization problem with penalty in
space $V_{0}$:
\begin{equation} \label{eq_PII_25}
F_{\theta } ({\bf u})=\frac{1}{2} A\, ({\bf u}, {\bf u})-L({\bf
u})+J_{\theta } ({\bf u})\to \mathop{\min }\limits_{ {\bf u} \,
\in V_{0} } .
\end{equation}
Note, that the introduction of the penalty corresponds to the
introduction of a conditional intermediate Winkler layer between
the bodies with the stiffness coefficient ${1\mathord{\left/
{\vphantom {1 \theta }} \right. \kern-\nulldelimiterspace} \theta
}$. The quantity $\sigma _{\alpha \beta n} =\sigma _{\alpha \, n}
=\sigma _{\beta \, n} ={(d_{\alpha \beta } -u_{\alpha \, n}
-u_{\beta \, n} )^{-} \mathord{\left/ {\vphantom {(d_{\alpha \beta
} -u_{\alpha \, n} -u_{\beta \, n} )^{-}  \theta }} \right.
\kern-\nulldelimiterspace} \theta }$ has a sense of the normal
contact stress between the bodies $\Omega_{\alpha}$ and
$\Omega_{\beta}$, and the penalty $J_{\theta} ({\bf u})$
represents the total work of the normal contact stress.

Now consider the properties of the penalty term (\ref {eq_PII_24})
in more detail. The functional $J_{\theta} ({\bf u})$ is
nonnegative
\begin{equation} \label{eq_PII_26}
\left(\forall {\bf u} \in V_{0} \right)\left\{J_{\theta } ({\bf
u})\ge 0\right\},
\end{equation}
and G\^{a}teaux differentiable in $V_{0}$:
\begin{equation} \label{eq_PII_27}
J'_{\theta } ({\bf u}, {\bf v})=-\frac{1}{\theta } \sum
_{\left\{\alpha ,\, \, \beta \right\}\, \in \, Q} \int _{S_{\alpha
\beta} }\left(d_{\alpha \beta }-u_{\alpha \, n}-u_{\beta \, n}
\right)^{-}\, (v_{\alpha \, n} +v_{\beta \, n} )\, dS ,
\end{equation}
Moreover, the G\^{a}teaux differential $J'_{\theta} ({\bf u}, {\bf
v})$ is linear in $\bf v$ and nonlinear in $\bf u$.

\textbf{Lemma~2.} \textit{If the surfaces $S_{\alpha
\beta}$},\textit{ $\left\{\alpha ,\beta \right\}\in Q,$ are
Lipschitz and $d_{\alpha \beta} \in H_{00}^{1/2}(\Xi_{\alpha
})$},\textit{ then $J'_{\theta} ({\bf u},{\bf v})$ satisfies the
following properties:}
\begin{equation} \label{eq_PII_29}
\left(\forall {\bf u} \in V_{0} \right)\left(\exists
\tilde{R}>0\right)\left(\forall {\bf v}\in V_{0} \right)\left\{\,
\left|J'_{\theta } ({\bf u}, {\bf v})\right|\le \tilde{R} \left\|
{\bf v}\right\| _{V_{0}} \right\},
\end{equation}
\begin{equation} \label{eq_PII_30}
\left(\forall {\bf u}, {\bf v}\in V_{0} \right)\left\{\,
J'_{\theta } ({\bf u}+{\bf v},{\bf v})-J'_{\theta } ({\bf u}, {\bf
v})\ge 0\right\},
\end{equation}
\begin{equation} \label{eq_PII_31}
\left(\exists D>0\right)\left(\forall {\bf u}, {\bf v},{\bf w} \in
V_{0} \right)\left\{\, \left|J'_{\theta } ({\bf u}+{\bf w},{\bf
v})-J'_{\theta } ({\bf u}, {\bf v})\right|\le D \left\| {\bf
v}\right\| _{V_{0}} \left\| {\bf w}\right\| _{V_{0}} \right\}.
\end{equation}

\textbf{Proof.} At first, let us show the satisfaction of property
(\ref {eq_PII_29}). Let us write $J'_{\theta } ({\bf u}, {\bf v})$
in the extended form: $J'_{\theta } ({\bf u}, {\bf v})=\sum
_{\left\{\alpha ,\, \, \beta \right\}\, \in \, Q} j'_{\alpha \beta
} ({\bf u}, {\bf v})$, where
\begin{equation} \label{eq_PII_32}
j'_{\alpha \beta } ({\bf u}, {\bf v})=-\frac{1}{\theta }
\left(\int _{S_{\alpha \beta } }g_{\alpha \beta } ({\bf u})\, \,
{\bf n}_{\alpha } \cdot {\rm{Tr}}_{\alpha }^{0} ({\bf v}_{\alpha})
\, dS +\int _{S_{\alpha \beta } }g_{\alpha \beta } ({\bf u})\, \,
{\bf n}_{\beta } \cdot  {\rm{Tr}}_{\beta}^{0} ({\bf v}_{\beta}) \,
dS \right),
\end{equation}
and $g_{\alpha \beta } ({\bf u})=\left(d_{\alpha \beta }-u_{\alpha
\, n}-u_{\beta \, n} \right)^{-}$.

Taking into account the following inequality for real numbers:
\begin{equation} \label{eq_PII_33}
\left(\sum _{i=1}^{m}c_{i}  \right)^{2} \le m\sum
_{i=1}^{m}c_{i}^{2} , \, \,\, c_{i} \in {\mathbb R}, \,\,\,
i=1,2,...,m , \,\,\, m\in {\mathbb N},
\end{equation}
and the Schwarz inequality, we obtain
\begin{equation} \label{eq_PII_34}
\left(\int _{S_{\alpha \beta } }g_{\alpha \beta } ({\bf u})\, \,
{\bf n}_{\alpha } \cdot {\rm{Tr}}_{\alpha }^{0} ({\bf v}_{\alpha})
\, dS \right)^{2} \le 3 \, q_{\alpha \beta } ({\bf u})\left\| \,
{\rm{Tr}}_{\alpha }^{0} ({\bf v}_{\alpha}) \right\| _{[L_{2}
(S_{\alpha \beta } )]^{m} }^{2},
\end{equation}
where $q_{\alpha \beta } ({\bf u})=\left\| g_{\alpha \beta } ({\bf
u})\, \, {\bf n}_{\alpha } \right\| _{[L_{2} (S_{\alpha \beta }
)]^{3} }^{2} +\varepsilon _{\alpha \beta } $, $\varepsilon
_{\alpha \beta } >0$, $\beta \in B_{\alpha } $, $\alpha
=1,2,...,N$.

From the trace theorems in \cite {PII_22}, it follows, that
\begin{equation} \label{eq_PII_35}
\left(\exists T_{\alpha \beta } >0\right)\left(\forall {\bf
v}_{\alpha } \in V_{\alpha }^{0} \right)\, \left\{\, \left\| {\rm
Tr}_{\alpha }^{0} ({\bf v}_{\alpha } )\right\| _{[L_{2} (S_{\alpha
\beta } )]^{3} }^{2} \le T_{\alpha \beta }^{2} \left\| {\bf
v}_{\alpha } \right\| _{V_{\alpha } }^{2} \right\}.
\end{equation}
Substituting (\ref {eq_PII_35}) into (\ref {eq_PII_34}), we come
to an inequality $\left|\int _{S_{\alpha \beta } }g_{\alpha \beta
} ({\bf u})\, \, {\bf n}_{\alpha } \cdot {\rm {Tr}}_{\alpha }^{0}
({\bf v}_{\alpha})\, dS \right|\le \tilde{s}_{\alpha \beta } ({\bf
u})\left\| {\bf v}_{\alpha } \right\| _{V_{\alpha } } $, where
$\tilde{s}_{\alpha \beta } ({\bf u})=T_{\alpha \beta }
\sqrt{3q_{\alpha \beta } ({\bf u})}
>0$, $\beta \in B_{\alpha } $, $\alpha =1,2,...,N$.

Similarly to this, we obtain the same inequality for the second
term of relationship (\ref {eq_PII_32}). Hence
\[\left|j'_{\alpha \beta } ({\bf u}, {\bf v})\right|\le \frac{1}{\theta }
\left(\tilde{s}_{\alpha \beta } ({\bf u})\left\| {\bf v}_{\alpha }
\right\| _{V_{\alpha } } +\tilde{s}_{\beta \alpha } ({\bf
u})\left\| {\bf v}_{\beta }
 \right\| _{V_{\beta } } \right).\]
As a result, we find
\[\left|J'_{\theta } ({\bf u}, {\bf v})\right|\le \sum _{\left\{\alpha ,\,
\, \beta \right\}\, \in \, Q}\left|j'_{\alpha \beta } ({\bf u},
{\bf v})\right| \le \tilde{R}({\bf u})\left\|{\bf v}\right\|
_{V_{0} } ,\] where $\tilde{R}({\bf u})=\frac{1}{\theta } \sum
_{\left\{\alpha ,\, \, \beta \right\}\, \in \,
Q}\left(\tilde{s}_{\alpha \beta } ({\bf u})+\tilde{s}_{\beta
\alpha } ({\bf u})\right)
>0$. Inequality (\ref {eq_PII_29}) is proved.

Now, we prove that condition (\ref {eq_PII_30}) holds. For this we
use the next inequality
\begin{equation} \label{eq_PII_36}
\left(\forall y,z\in {\mathbb R}\right)\left\{\, \left[(y-z)^{-}
-y^{-} \right]z\le 0\right\}.
\end{equation}
Rewrite $J'_{\theta } ({\bf u}+{\bf w},{\bf v})-J'_{\theta } ({\bf
u}, {\bf v})$ in the following way:
\begin{equation} \label{eq_PII_37}
J'_{\theta } ({\bf u}+{\bf w},{\bf v})-J'_{\theta } ({\bf u}, {\bf
v})=-\frac{1}{\theta } \sum _{\left\{\alpha ,\, \, \beta
\right\}\,
 \in \, Q}h_{\alpha \beta } ({\bf u},{\bf w},{\bf v}),
\end{equation}
where $h_{\alpha \beta } ({\bf u},{\bf w},{\bf v})=\int
_{S_{\alpha \beta } }r_{\alpha \beta } ({\bf u},{\bf w})\,
\left(v_{\alpha \, n} +v_{\beta \, n} \right)\, dS $, $r_{\alpha
\beta } ({\bf u},{\bf w})=g_{\alpha \beta } ({\bf u}+{\bf
w})-g_{\alpha \beta } ({\bf u}),\, {\bf x} \in S_{\alpha \beta }.$

In view of property (\ref {eq_PII_36}), we obtain
\[\left(\forall {\bf u}, {\bf v}\in V_{0} \right)\left\{r_{\alpha \beta }
 ({\bf u}, {\bf v})\, \left(v_{\alpha \, n} +v_{\beta \, n} \right)\le 0,\,
  \, \, {\bf x} \in S_{\alpha \beta } \right\}.\]
Therefore, since $\theta >0$, we come to an inequality
\[J'_{\theta } ({\bf u}+{\bf v},{\bf v})-J'_{\theta }
 ({\bf u}, {\bf v})=-\frac{1}{\theta } \sum _{\left\{\alpha ,\,
 \, \beta \right\}\, \in \, Q}\int _{S_{\alpha \beta } }r_{\alpha \beta }
  ({\bf u}, {\bf v})\, \left(v_{\alpha \, n} +v_{\beta \, n} \right)\,
   dS \ge 0 ,\,\, \forall {\bf u}, {\bf v}\in V_{0} .\]

Finally, let us prove the satisfaction of condition (\ref
{eq_PII_31}). Using the next inequality for real numbers
\begin{equation} \label{eq_PII_38}
 \left(\forall y,z\in {\mathbb R}\right)\left\{\, \left|y^{-} -z^{-}
\right|\le \left|y-z\right|\right\},
\end{equation}
we obtain
\begin{equation} \label{eq_PII_39}
r_{\alpha \beta } ({\bf u},{\bf w})\le \left|w_{\alpha \, n}
+w_{\beta \, n} \right|, \,\, {\bf x} \in S_{\alpha \beta } .
\end{equation}
Taking into account, that the components of the outer unit normal
to $\partial \Omega_{\alpha}$ satisfy the property
$\mathop{\max}\limits_{j=1,2,3} \left| n_{\alpha \, j} \right| \le
1$, and using inequalities (\ref {eq_PII_33}) and (\ref
{eq_PII_39}), we find
\begin{equation} \label{eq_PII_40}
\int _{S_{\alpha \beta }} r_{\alpha \beta }^{2} ({\bf u},{\bf w})
\, dS \le 6 \left(\left\| {\rm{Tr}}_{\alpha }^{0} ({\bf
w}_{\alpha})\right\| _{[L_{2} (S_{\alpha \beta } )]^{3} }^{2}
+\left\| {\rm{Tr}}_{\beta }^{0}
 ({\bf w}_{\beta})\right\| _{[L_{2} (S_{\alpha \beta } )]^{3} }^{2}\right).
\end{equation}

Now, let us write $h_{\alpha \beta } ({\bf u},{\bf w},{\bf v})$ as
follows:
\[h_{\alpha \beta } ({\bf u},{\bf w},{\bf v})=\int _{S_{\alpha \beta } }r_{\alpha \beta }
({\bf u},{\bf w})\, v_{\alpha \, n} \, dS +\int _{S_{\alpha \beta
} }r_{\alpha \beta } ({\bf u},{\bf w})\, v_{\beta \, n} \, dS .\]
Consider the term $\int _{S_{\alpha \beta } }r_{\alpha \beta }
({\bf u},{\bf w})\, v_{\alpha \, n} \, dS $ in more detail. Let us
use inequalities (\ref {eq_PII_33}), (\ref {eq_PII_40}), and the
Schwarz inequality:
\[\left(\int _{S_{\alpha \beta } }r_{\alpha \beta } ({\bf u},{\bf w})\, v_{\alpha \, n} \,
 dS \right)^{2} \le \int _{S_{\alpha \beta } }r_{\alpha \beta }^{2} ({\bf u},{\bf w})\,
 dS \int _{S_{\alpha \beta } }v_{\alpha \, n}^{2} \, dS \le \]
\[\le 18\left(\, \left\| {\rm{Tr}}_{\alpha }^{0} ({\bf w}_{\alpha})\right\|
_{[L_{2} (S_{\alpha \beta } )]^{3} }^{2} +\left\| {\rm{Tr}}_{\beta
}^{0}
 ({\bf w}_{\beta})\right\| _{[L_{2} (S_{\alpha \beta } )]^{3} }^{2} \right)\,
 \, \left\| {\rm{Tr}}_{\alpha }^{0} ({\bf v}_{\alpha})\right\| _{[L_{2}
  (S_{\alpha \beta } )]^{3} }^{2} .\]
Using inequality (\ref {eq_PII_35}), we find further
\[ \left|\int _{S_{\alpha \beta } }r_{\alpha \beta } ({\bf u},{\bf w})\, v_{\alpha \, n} \, dS
\right|\le T_{\alpha \beta 1} \left\| {\bf v}_{\alpha } \right\|
_{V_{\alpha } } \left(\, \left\| {\bf w}_{\alpha } \right\|
_{V_{\alpha } } +\left\| {\bf w}_{\beta } \right\| _{V_{\beta } }
\right), \] where $T_{\alpha \beta 1} =3\, T_{\alpha \beta } \,
T_{\alpha \beta }^{*} \sqrt{2} \,
>0$, $T_{\alpha \beta }^{*} =\max \left\{T_{\alpha \beta } ,T_{\beta \alpha }
\right\}>0$. In much the same way we come to an inequality
\[\left|\int _{S_{\alpha \beta } }r_{\alpha \beta } ({\bf u},{\bf w})\,
 v_{\beta \, n} \, dS \right|\le T_{\alpha \beta 2} \left\| {\bf v}_{\beta }
  \right\| _{V_{\beta } } \left(\, \left\| {\bf w}_{\alpha } \right\| _{V_{\alpha } }
  +\left\| {\bf w}_{\beta } \right\| _{V_{\beta } } \right),\]
where $T_{\alpha \beta 2} =3\,T_{\beta \alpha } \, T_{\alpha \beta
}^{*} \sqrt{2} >0$.

Taking into account the last two inequalities, we establish
\[\left|h_{\alpha \beta } ({\bf u},{\bf w},{\bf v})\right|\le \left|\int _{S_{\alpha \beta } }
r_{\alpha \beta } ({\bf u},{\bf w})\, v_{\alpha \, n} \, dS
\right|+ \left|\int _{S_{\alpha \beta } }r_{\alpha \beta } ({\bf
u},{\bf w})\, v_{\beta \, n} \, dS \right|\le \]
\[ \le C_{\alpha \beta } \left(\left\| {\bf v}_{\alpha } \right\| _{V_{\alpha } } \left\|
{\bf w}_{\alpha } \right\| _{V_{\alpha } } +\left\| {\bf
v}_{\alpha } \right\| _{V_{\alpha } } \left\| {\bf w}_{\beta }
\right\| _{V_{\beta } } +\left\| {\bf v}_{\beta } \right\|
_{V_{\beta } } \left\| {\bf w}_{\alpha } \right\| _{V_{\alpha } }
+\left\| {\bf v}_{\beta } \right\| _{V_{\beta } } \left\| {\bf
w}_{\beta } \right\| _{V_{\beta } } \right), \] where $C_{\alpha
\beta } =\max \left\{T_{\alpha \beta 1} ,T_{\alpha \beta 2}
\right\}>0$.

As a result, we obtain
\[\left|J'_{\theta } ({\bf u}+{\bf w},{\bf v})-J'_{\theta }
({\bf u}, {\bf v})\right|\le \frac{1}{\theta } \sum _{\alpha ,\,
\beta =1}^{N}\left|h_{\alpha \beta }
 ({\bf u},{\bf w},{\bf v})\right| \le \]
\[\le \frac{2C}{\theta } \left(N\sum _{\alpha =1}^{N}\left\| {\bf v}_{\alpha } \right\|
_{V_{\alpha } } \left\| {\bf w}_{\alpha } \right\| _{V_{\alpha } }
 +\sum _{\alpha =1}^{N}\left\| {\bf v}_{\alpha } \right\| _{V_{\alpha } }
 \sum _{\beta =1}^{N}\left\| {\bf w}_{\beta } \right\| _{V_{\beta } }  \right)\le \]
\[\le \frac{2C\left(N+1\right)}{\theta }
\sum _{\alpha =1}^{N}\left\| {\bf v}_{\alpha } \right\|
_{V_{\alpha } }
 \sum _{\beta =1}^{N}\left\| {\bf w}_{\beta } \right\| _{V_{\beta } }
 \le D\left\|{\bf v}\right\| \, \, \left\| {\bf w}\right\| , \,\,\, \forall {\bf u}, {\bf v},{\bf w} \in V_{0} ,\]
where $D={2C\left(N+1\right)\mathord{\left/ {\vphantom
{2C\left(N+1\right) \theta }} \right. \kern-\nulldelimiterspace}
\theta } >0$, $C=\mathop{\max }\limits_{1\le \alpha
,\beta \le N} C_{\alpha \beta } >0$. $\Box$ \\

\textbf{Theorem~2.} \textit{Suppose that the conditions of
\textbf{\textit{Lemma~1}} and \textbf{\textit{Lemma~2}} hold. Then
there exists a unique solution of the nonquadratic minimization
problem (\ref {eq_PII_25}) in $V_{0},$ and this problem is
equivalent to the following nonlinear variational equation:}
\begin{equation} \label{eq_PII_41}
F'_{\theta } ({\bf u}, {\bf v})=A\, ({\bf u}, {\bf v})+J'_{\theta
} ({\bf u}, {\bf v})-L\, ({\bf v})=0, \,\,\, \forall \, {\bf v}
\in V_{0}, \,\,\, {\bf u} \in V_{0} .
\end{equation}

\textbf{Proof.} As shown in \cite {PII_7}, due to properties (\ref
{eq_PII_18})~--~(\ref {eq_PII_21}), the functional $F({\bf u})$ is
strictly convex and coercive ($\mathop{\lim }\limits_{\left\|{\bf
u}\right\| _{V_{0} } \to \infty } F({\bf u})=\infty $), and the
differential $F'({\bf u}, {\bf v})$ is linear and continuous in
$\bf v$.

From property (\ref {eq_PII_30}), it follows that the penalty term
$J_{\theta } ({\bf u})$ is convex in $V_{0} $ \cite {PII_15}.

Now consider the properties of the functional
\[F_{\theta } ({\bf u})=F({\bf u})+J_{\theta } ({\bf u}), \,\, {\bf u} \in V_{0} .\]
This functional is G\^{a}teaux differentiable in $V_{0}$:
\[F'_{\theta } ({\bf u}, {\bf v})=F'({\bf u}, {\bf v})+J'_{\theta } ({\bf u}, {\bf v}), \,\, {\bf u},{\bf v} \in V_{0} ,\]
and strictly convex, as the sum of the convex functional
$J_{\theta} ({\bf u})$ and the strictly convex functional $F({\bf
u})$. In addition, since the functionals $F'({\bf u},{\bf v})$ and
$J'_{\theta} ({\bf u},{\bf v})$ are linear and continuous in $\bf
v$, it follows that $F'_{\theta } ({\bf u}, {\bf v})$ is also
linear and continuous in $\bf v$. Hence, according to \cite
{PII_4}, the functional $F_{\theta} ({\bf u})$ is weakly lower
semicontinuous. Due to the coercivity of $F({\bf u})$ and property
(\ref {eq_PII_26}), we obtain that $\mathop{\lim
}\limits_{\left\|{\bf u}\right\| _{V_{0} } \to \infty } F_{\theta
} ({\bf u})=\infty$.

Since the functional $F_{\theta}({\bf u})$ is weakly lower
semicontinuous, coercive, strictly convex, and G\^{a}teaux
differentiable in the closed reflexive Banach space $V_{0}$, then
according to theorems in \cite {PII_4}, there exists a unique
solution of the minimization problem (\ref {eq_PII_25}) in
$V_{0}$, and this problem is equivalent to the variational
equation (\ref {eq_PII_41}). $\Box$

Now let us prove that the solution of the penalty variational
equation (\ref {eq_PII_41}) converges strongly to the solution of
the original variational inequality (\ref {eq_PII_23}) as $\theta
\to 0$.

Let us rewrite the variational equation (\ref {eq_PII_41}) in the
following equivalent form
\begin{equation} \label{eq_PII_42}
F'_{\theta } ({\bf u}, {\bf v})=A\, ({\bf u}, {\bf v})-\langle
L,{\bf v}\rangle +\frac{1}{\theta } \left\langle \Phi ({\bf
u}),{\bf v}\right\rangle =0, \,\, \forall \, {\bf v}\in V_{0},
\,\, {\bf u}\in V_{0} ,
\end{equation}
where $\langle Y,{\bf u}\rangle =Y({\bf u})$ is the action of a
functional $Y\in V_{0}^{*}$ on an element ${\bf u}\in V_{0}$,
$V_{0}^{*}$ is the space dual to $V_{0}$, $\Phi =\Psi':\, \, V_{0}
\to V_{0}^{*}$ is the G\^{a}teaux derivative of the functional
$\Psi ({\bf u})=\theta J_{\theta} ({\bf u})$, and $\left\langle
\Phi ({\bf u}),{\bf v}\right\rangle =\left\langle \Psi'({\bf
u}),{\bf v}\right\rangle =\theta J'_{\theta} ({\bf u}, {\bf v})$.

We have proved the next lemma.

\textbf{Lemma~3.} \textit{Suppose that the conditions of
\textbf{\textit{Lemma~2}} hold. Then the operator $\Phi: \,\,
V_{0} \to V_{0}^{*}$ in problem (\ref {eq_PII_42}) is a penalty
operator for the kinematically allowable displacements set $K$,
i.e.}

1).\textit{~$\Phi $ is monotone in $V_{0} $:}

\[\left(\forall {\bf u}, {\bf v}\in V_{0} \right)\left\{\left\langle
\Phi ({\bf u})-\Phi (\bf{v}),u-{\bf v}\right\rangle \ge
0\right\};\]

2).\textit{~$\Phi $ satisfies the Lipschitz condition in $V_{0}
$:}

\[\left(\exists C>0\right)\left(\forall {\bf u}, {\bf v}\in V_{0}
 \right)\left\{\, \left\| \Phi ({\bf u})-\Phi (\bf{v})\right\| _{V_{0}^{*} }
\le C\left\| {\bf u}-{\bf v}\right\|_{V_{0}} \right\};\]

3).\textit{~The kernel of operator $\Phi $ is equal to the set
$K$:}

\[{\rm Ker}\, (\Phi )=\left\{ {\bf u}:\, \, \, \, {\bf u}\in V_{0} ,\,
\, \, \Phi ({\bf u})=0\right\}=K.\]

\textbf{Proof.} The monotonicity of operator $\Phi$ follows from
condition (\ref {eq_PII_30}) and the satisfaction of Lipschitz
condition follows from property (\ref {eq_PII_31}).

If ${\bf u}\in K$, then $\Phi ({\bf u})\equiv 0$, and, on the
contrary, if $\Phi ({\bf u})\equiv 0$, we have ${\bf u}\in K$.
Hence, ${\rm Ker}\, (\Phi )=K$. For more details see \cite
{PII_7}. $\Box$

Now, using the results of works \cite {PII_21,PII_11,PII_14}, let
us prove the proposition on the strong convergence of the penalty
method, applied to the variational inequality (\ref {eq_PII_23}).

\textbf{Theorem~3.} \textit{Suppose that the conditions of
\textbf{\textit{Lemma~1}} and \textbf{\textit{Lemma~2}} hold,}
$\bar{\bf u}\in K$ \textit{is a unique solution of the variational
inequality (\ref {eq_PII_23}), and $\bar{\bf u}_{\, \theta } \in
V_{0} $ is a unique solution of the penalty variational equation
(\ref {eq_PII_42}) with the penalty parameter $\theta
>0$. Then $\bar{\bf u}_{\, \theta } \mathop{\to }\limits_{\theta
\to 0} \bar{\bf u}$ strongly in $V_{0} $},\textit{ i.e. $\left\|
\bar{\bf u}_{\, \theta } -\bar{\bf u}\right\|_{V_{0}}  \mathop{\to
}\limits_{\theta \to 0} 0$}.

\textbf{Proof.} In works \cite {PII_21,PII_11} it is proved that
if conditions (\ref {eq_PII_19})~--~(\ref {eq_PII_21}) hold, $\Phi
$ is a penalty operator for the set $K$, and there exist the
solutions of problems (\ref {eq_PII_23}) and (\ref {eq_PII_42}),
then the sequence $\{ \bar{\bf u}_{\, \theta } \}$ is bounded:

\[\left(\exists \, \tilde{C}\in (0;\infty )\right)\left(\forall \theta >0\right)\left\{\,
\left\| \bar{\bf u}_{\, \theta } \right\| _{V_{0} } \le
\tilde{C}\right\},\] and there exists such subsequence $\{
\bar{\bf u}_{\,\theta _{1} } \} \subset \{ \bar{\bf u}_{\, \theta
} \}$, which converges weakly in $V_{0}$ to some solution of
variational inequality (\ref {eq_PII_23}), i.e.
\[\left(\exists \, \{ \bar{\bf u}_{\,\theta _{1} } \}
\subset \{ \bar{\bf u}_{\, \theta } \} \right)\left(\forall Y\in
V_{0}^{*}
 \right)\left\{\, \left\langle Y,\bar{\bf u}_{\,\theta _{1} }
  \right\rangle \mathop{\to }\limits_{\theta _{1} \to 0} \left\langle
   Y,\bar{\bf u}\right\rangle \right\}.\]

Moreover, in \cite {PII_21,PII_11} it is shown that any weakly
convergent subsequence of the sequence $\{\bar{\bf u}_{\,
\theta}\}$ converges weakly in $V_{0} $ to some solution of
variational problem (\ref {eq_PII_23}).

Now let us assume that variational problems (\ref {eq_PII_23}) and
(\ref {eq_PII_42}) have unique solutions. Then, as follows from
above, the sequence $\{ \bar{\bf u}_{\, \theta } \} $ has a unique
partial weak limit $\bar{\bf u}\in K$.

Since the sequence $\{ \bar{\bf u}_{\, \theta } \} $ has a unique
weak limit point, and is bounded, then according to the theorem in
\cite {PII_4}, it is weakly convergent to this point, i.e.
\begin{equation} \label{eq_PII_43}
\left(\forall \, Y\in V_{0}^{*} \right)\left\{\, \left\langle
Y,\bar{\bf u}_{\, \theta } \right\rangle \mathop{\to
}\limits_{\theta \to 0} \left\langle Y,\bar{\bf u}\right\rangle
\right\},
\end{equation}
where $\bar{\bf u} \in K$ is a unique solution of variational
inequality (\ref {eq_PII_23}).

Further, let us show that $\{ \bar{\bf u}_{\, \theta } \}$
converges strongly to $\bar{\bf u}\in K$ as $\theta \to 0$.

Due to (\ref {eq_PII_43}), we get
\begin{equation} \label{eq_PII_44}
\left\langle L,\bar{\bf u}_{\, \theta } -\bar{\bf u}\right\rangle
\mathop{\to } \limits_{\theta \to 0} 0, \, \, A ({\bf v},\bar{\bf
u}_{\, \theta } -\bar{\bf u})=\left\langle A'_{1}
(\bf{v}),\bar{\bf u}_{\, \theta } -\bar{\bf u}\right\rangle
\mathop{\to }\limits_{\theta \to 0} 0, \,\, \forall \, {\bf v}\in
V_{0} ,
\end{equation}
where $A'_{1} (\bf{v})$ is the G\^{a}teaux derivative of the
functional $A_{1} ({\bf v})=\frac{1}{2} A\, ({\bf v},{\bf v})$,
${\bf v}\in V_{0}$.

Since $\bar{\bf u}_{\, \theta } $ is a solution of the penalty
variational equation (\ref {eq_PII_42}), it is obvious that
\[\left\langle L,\bar{\bf u}_{\, \theta } -{\bf v}\right\rangle =A\,
(\bar{\bf u}_{\, \theta } ,\bar{\bf u}_{\, \theta } -{\bf
v})+\frac{1}{\theta }
 \left\langle \Phi (\bar{\bf u}_{\, \theta } ),\bar{\bf u}_{\, \theta } -{\bf v}\right\rangle, \,\,
 \forall \, {\bf v}\in K.\]

Taking into account the monotonicity of the penalty operator
$\Phi$ and the property $\left(\forall \, {\bf v}\in
K\right)\left\{\Phi ({\bf v})=0\right\}$, we obtain
\[\left\langle L,\bar{\bf u}_{\, \theta } -{\bf v}\right\rangle =A\, (\bar{\bf u}_{\, \theta }
 ,\bar{\bf u}_{\, \theta } -{\bf v})+\frac{1}{\theta } \left\langle \Phi (\bar{\bf u}_{\,
  \theta } )-\Phi (\bf{v}),\bar{\bf u}_{\, \theta }
  -{\bf v}\right\rangle \ge A\, (\bar{\bf u}_{\, \theta }
  ,\bar{\bf u}_{\, \theta } -{\bf v}), \,\, \forall \, {\bf v}\in K.\]
In view of this property and the nonnegativity of the bilinear
form $A$, we get an inequality \noindent $A\, (\bar{\bf
u},\bar{\bf u}_{\, \theta } -\bar{\bf u})-\left\langle L,\bar{\bf
u}_{\, \theta } -\bar{\bf u}\right\rangle \le A\, (\bar{\bf u}_{\,
\theta } ,\bar{\bf u}_{\, \theta } -\bar{\bf u})-\left\langle
L,\bar{\bf u}_{\, \theta } -\bar{\bf u}\right\rangle \le 0$. Hence
\begin{equation} \label{eq_PII_45}
A\, (\bar{\bf u},\bar{\bf u}_{\, \theta } -\bar{\bf u})\le A\,
(\bar{\bf u}_{\, \theta } ,\bar{\bf u}_{\, \theta } -\bar{\bf
u})\le \left\langle L,\bar{\bf u}_{\, \theta } -\bar{\bf
u}\right\rangle .
\end{equation}

Passing to the limit in expression (\ref {eq_PII_45}) as $\theta
\to 0$, and taking into account property (\ref {eq_PII_44}), we
obtain
\[A\, (\bar{\bf u}_{\, \theta } ,\bar{\bf u}_{\, \theta } -\bar{\bf u})\mathop{\to }\limits_{\theta \to 0} 0.\]

Further, in view of the coercivity of bilinear form $A$, it
follows that
\[0\le B\left\| \bar{\bf u}_{\, \theta } -\bar{\bf u}\right\| _{V_{0} }^{2} \le A\,
(\bar{\bf u}_{\, \theta } ,\bar{\bf u}_{\, \theta } -\bar{\bf
u})-A\, (\bar{\bf u},\bar{\bf u}_{\, \theta } -\bar{\bf
u})\mathop{\to }\limits_{\theta \to 0} 0, \,\, B>0.\] As a result,
we establish that $\left\| \bar{\bf u}_{\, \theta } -\bar{\bf
u}\right\| _{V_{0} } \mathop{\to }\limits_{\theta \to 0} 0$.
$\Box$

Thus, using the penalty method we reduced the solution of the
original variational inequality (\ref {eq_PII_23}) on the closed
convex set $K$ to the solution of the nonlinear variational
equation (\ref {eq_PII_41}) in the whole space $V_{0}$, which
depends on the penalty parameter $\theta >0$. We also proved the
existence of a unique solution of the penalty variational equation
(\ref {eq_PII_41}) and its strong convergence to a solution of the
original variational inequality (\ref {eq_PII_23}) as the penalty
parameter $\theta$ tends to zero.

In the following section let us consider some iterative methods to
solve such nonlinear variational equations.

\section{Iterative methods for nonlinear variation equations}

Consider an abstract nonquadratic minimization problem in form
(\ref {eq_PII_25}), and equivalent nonlinear in $\bf u$
variational equation in form (\ref {eq_PII_41}), where $V_{0}$ is
some closed reflexive Banach space, $A\,({\bf u},{\bf v})$ is a
bilinear form in $V_{0}$, $L\,({\bf v})$ is a linear functional,
and the term $J'_{\theta}({\bf u},{\bf v})$ is linear in $\bf v$
and nonlinear in $\bf u$. Suppose that conditions (\ref
{eq_PII_18})~--~(\ref {eq_PII_21}), (\ref {eq_PII_29})~--~(\ref
{eq_PII_31}) are satisfied. Hence, there exists a unique solution
of problem (\ref {eq_PII_41}).

For the numerical solution of the nonlinear variational equation
(\ref {eq_PII_41}) let us use the following iterative method \cite
{PII_23,PII_7,PII_8,PII_47}:
\begin{equation} \label{eq_PII_46}
G\, ({\bf u}^{k+1} ,{\bf v})=G\, ({\bf u}^{k} ,{\bf v})-\gamma
\left[A\, ({\bf u}^{k} ,{\bf v})+J'_{\theta } ({\bf u}^{k} ,{\bf
v})-L\, (\bf{v})\right], \,\, k=0,1,...\,,
\end{equation}
where $G\,({\bf u},{\bf v})$ is some bilinear form assigned in
$V_{0}$, ${\bf u}^{k} \in V_{0}$, $k=1,2,...$ is the $k$-th
approximation to the exact solution $\bar{\bf u} \in V_{0}$ of
problem (\ref {eq_PII_41}), ${\bf u}^{0} \in V_{0}$ is an initial
approximation, and $\gamma \in {\mathbb R}$ is an iterative
parameter.

This iterative method can be viewed as a descent method for the
minimization problem (\ref {eq_PII_25}) with the choice of the
descent direction via an auxiliary operator \cite{PII_4}. On the
other hand, this method can be viewed as an implicit successive
iteration method for the variational equation (\ref {eq_PII_41}),
and the bilinear form $G\,({\bf u},{\bf v})$ can be interpreted as
a preconditioner.

Using the methodology, developed in \cite{PII_11} for the case of
linear variational equations, we have proved the following theorem
on the convergence of the iterative method (\ref {eq_PII_46}) for
nonlinear variational equations.

\textbf{Theorem~4.} \textit{Suppose that the bilinear form
$G\,({\bf u},{\bf v})$ is symmetric, continuous and coercive:}
\begin{equation} \label{eq_PII_47}
\left(\forall {\bf u}, {\bf v}\in V_{0} \right)\left\{G\, ({\bf
u}, {\bf v})=G\, ({\bf v}, {\bf u})\right\},
\end{equation}
\begin{equation} \label{eq_PII_48}
\left(\exists \tilde{M}>0\right)\left(\forall {\bf u}, {\bf v}\in
V_{0} \right)\left\{\, \left|G\, ({\bf u}, {\bf v})\right|\le
\tilde{M}\left\|{\bf u}\right\| _{V_{0} } \left\|{\bf v}\right\|
_{V_{0} } \right\},
\end{equation}
\begin{equation} \label{eq_PII_49}
\left(\exists \tilde{B}>0\right)\left(\forall {\bf u} \in V_{0}
\right)\left\{G\, ({\bf u}, {\bf u})\ge \tilde{B}\left\|{\bf
u}\right\| _{V_{0} }^{2} \right\},
\end{equation}
\textit{properties (\ref {eq_PII_19})~--~(\ref {eq_PII_21}), (\ref
{eq_PII_29})~--~(\ref {eq_PII_31}) are satisfied, and the
iterative parameter lies in the interval $\gamma \in (0;\gamma
_{2} )$}, $\gamma _{2} ={2B\tilde{B}\mathord{\left/ {\vphantom
{2B\tilde{B} M_{*}^{2} }} \right. \kern-\nulldelimiterspace}
M_{*}^{2} } $, $M_{*} =M+D$.\textit{} \textit{Then the sequence
$\{ {\bf u}^{k} \}$}, \textit{obtained by the iterative method
(\ref {eq_PII_46}), converges strongly in $V_{0}$ to the exact
solution $\bar{\bf u} \in V_{0}$ of the variational equation (\ref
{eq_PII_41}), i.e. $\left\| {\bf u}^{k} -\bar{\bf u}\right\|
_{V_{0} } \mathop{\to }\limits_{k\to \infty } 0$},\textit{ and the
convergence rate  in the energy norm $\left\|{\bf u}\right\| _{G}
=\sqrt{G\, ({\bf u}, {\bf u})} $ is linear:}
\begin{equation} \label{eq_PII_50}
\left\| {\bf u}^{k+1} -\bar{\bf u}\right\| _{G} \le q\left\| {\bf
u}^{k} -\bar{\bf u}\right\| _{G},\,\, q=\sqrt{1-{\gamma
\left(2B-\gamma M_{*}^{2} /\tilde{B}\right)\mathord{\left/
{\vphantom {\gamma \left(2B-\gamma M_{*}^{2} /\tilde{B}\right)
\tilde{M}}} \right. \kern-\nulldelimiterspace} \tilde{M}} \, } <1.
\end{equation}
\textit{Moreover, the maximal convergence rate reaches as $\gamma
=\bar{\gamma }={B\tilde{B}\mathord{\left/ {\vphantom {B\tilde{B}
M_{*}^{2} }} \right. \kern-\nulldelimiterspace} M_{*}^{2} }
$}.\textit{}

\textbf{Proof.} Since the bilinear form $G\,({\bf u},{\bf v})$ is
symmetric, continuous and coercive, we may introduce a scalar
product and a norm
\[\left({\bf u},{\bf v}\right)_{G} =G\, ({\bf u}, {\bf v}),\,\,\,
  \left\|{\bf u}\right\| _{G} =\sqrt{G\, ({\bf u}, {\bf u})} , \,\,\, {\bf u},{\bf v}\in V_{0} .\]

From properties (\ref {eq_PII_48}) and (\ref {eq_PII_49}), it
follows that the norms $\left\| \, \cdot \, \right\| _{G}$ and
$\left\| \, \cdot \, \right\| _{V_{0}}$ are equivalent in space
$V_{0} $.

In each step $k\in \left\{0,1,...\right\}$ of method (\ref
{eq_PII_46}), we have to solve the linear variational problem:
\begin{equation} \label{eq_PII_51}
G\, ({\bf u}, {\bf v})=Y^{k} ({\bf v}), \,\,\, \forall \, {\bf
v}\in V_{0}, \,\,\, {\bf u}\in V_{0} ,
\end{equation}
where $Y^{k} ({\bf v})=G\, ({\bf u}^{k} ,{\bf v})-\gamma \left[A\,
({\bf u}^{k} ,{\bf v})+J'_{\theta } ({\bf u}^{k} ,{\bf v})-L\,
(\bf{v})\right]$ is linear in ${\bf v}$, ${\bf u}^{k} \in V_{0}$.
Using properties (\ref {eq_PII_19}), (\ref {eq_PII_21}), (\ref
{eq_PII_29}), and (\ref {eq_PII_48}), we obtain that $Y^{k}
(\bf{v})$ is continuous:
\begin{equation} \label{eq_PII_52}
\left(\exists {Z}_{k} >0\right)\left(\forall \, {\bf v}\in V_{0}
\right)\left\{\, \left|Y^{k} (\bf{v})\right|\le {Z}_{k}
\left\|{\bf v}\right\| _{V_{0} } \right\},
\end{equation}
where
 ${Z}_{k} =\tilde{M}\left\| {\bf u}^{k} \right\| _{V_{0} } +\left|\gamma
\right|\left(M\left\| {\bf u}^{k} \right\| _{V_{0} }
+\tilde{R}\,({\bf u}^{k} )+T\right)+\varepsilon
>0$, $\varepsilon >0$.

Since conditions (\ref {eq_PII_47}) -- (\ref {eq_PII_49}) and
(\ref {eq_PII_52}) are satisfied, we see that problem (\ref
{eq_PII_51}) has a unique solution ${\bf u}={\bf u}^{k+1} \in
V_{0} $.

Now let us show, that the sequence of solutions of problems (\ref
{eq_PII_51}) converges strongly to the solution of the original
variational equation (\ref {eq_PII_41}).

Suppose that $\bar{\bf u}\in V_{0} $ is the exact solution of
problem (\ref {eq_PII_41}). Introduce a notation ${\pmb \varphi
}^{k} :={\bf u}^{k} -\bar{\bf u}\in V_{0} $, $k=0,1,...$, and
rewrite (\ref {eq_PII_46}) as follows:
\[G\, (\bar{\bf u}+{\pmb \varphi }^{k+1} ,{\bf v})=G\, (\bar{\bf u}+{\pmb \varphi }^{k} ,
{\bf v})-\gamma \left[A\, (\bar{\bf u}+{\pmb \varphi }^{k} ,{\bf
v})+J'_{\theta } (\bar{\bf u}+{\pmb \varphi }^{k} ,{\bf v})-L\,
(\bf{v})\right].\]

Subtracting from this expression the identity $G\, (\bar{\bf
u},{\bf v}) \equiv G\, (\bar{\bf u},{\bf v})-\gamma [A\, (\bar{\bf
u},{\bf v})+J'_{\theta } (\bar{\bf u},{\bf v})-L\, (\bf{v})]$, we
obtain
\begin{equation} \label{eq_PII_53}
G\, ({\pmb \varphi }^{k+1} ,{\bf v})=G\, ({\pmb \varphi }^{k}
,{\bf v})-\gamma \left[A\, ({\pmb \varphi }^{k} ,{\bf
v})+J'_{\theta } (\bar{\bf u}+{\pmb \varphi }^{k} ,{\bf
v})-J'_{\theta } (\bar{\bf u},{\bf v})\right].
\end{equation}

Let us define the functional $H_{\theta } ({\bf u},{\bf w},{\bf
v})=J'_{\theta } ({\bf u}+{\bf w},{\bf v})-J'_{\theta } ({\bf u},
{\bf v})$, ${\bf u},{\bf v},{\bf w} \in V_{0} $, which is linear
in ${\bf v}$.

Due to properties (\ref {eq_PII_30}) and (\ref {eq_PII_31}), the
following conditions hold:
\begin{equation} \label{eq_PII_54}
\left(\forall {\bf u}, {\bf v}\in V_{0} \right)\left\{\, H_{\theta
} ({\bf u},{\bf v},{\bf v})\ge 0\right\},
\end{equation}
\begin{equation} \label{eq_PII_55}
\left(\exists D>0\right)\left(\forall {\bf u}, {\bf v},{\bf w} \in
V_{0} \right)\left\{\, \left|H_{\theta } ({\bf u},{\bf w},{\bf
v})\right|\le D\left\| {\bf w}\right\| _{V_{0} } \left\|{\bf
v}\right\| _{V_{0} } \right\}.
\end{equation}

Let us rewrite expression (\ref {eq_PII_53}) in the form
\begin{equation} \label{eq_PII_56}
G\, ({\pmb \varphi}^{k+1} -{\pmb \varphi}^{k} ,{\bf v})=-\gamma
\left[A\, ({\pmb \varphi }^{k} ,{\bf v})+H_{\theta } (\bar{\bf
u},\, {\pmb \varphi }^{k} ,\, {\bf v})\right], \,\, \forall \,
{\bf v}\in V_{0} .
\end{equation}

If we take ${\bf v}:={\pmb \varphi }^{k+1} -{\pmb \varphi }^{k}$
in (\ref {eq_PII_56}), we will have
\[\left\| {\pmb \varphi }^{k+1} -{\pmb \varphi }^{k} \right\| _{G}^{2}
 \le \left|\gamma \right|\left(\left|A\, ({\pmb \varphi }^{k} ,
 {\pmb \varphi }^{k+1} -{\pmb \varphi }^{k} )\right|+\left|H_{\theta }
  (\bar{\bf u},{\pmb \varphi }^{k} ,{\pmb \varphi }^{k+1} -{\pmb \varphi }^{k}
   )\right|\right).\]

Taking into account the continuity of bilinear form (\ref
{eq_PII_19}), and property (\ref {eq_PII_55}), we get
\[\left\| {\pmb \varphi }^{k+1} -{\pmb \varphi }^{k} \right\| _{G}^{2} \le \left|\gamma \right|M_{*} \left\| {\pmb \varphi }^{k} \right\| _{V_{0} } \left\| {\pmb \varphi }^{k+1} -{\pmb \varphi }^{k} \right\| _{V_{0} }, \,\, M_{*} =M+D>0.\]
Further, in view of the relation between norms
\begin{equation} \label{eq_PII_57}
\left\| {\pmb \varphi }^{k+1} -{\pmb \varphi }^{k} \right\|
_{V_{0} } \le {\left\| {\pmb \varphi }^{k+1} -{\pmb \varphi }^{k}
\right\| _{G} \mathord{\left/ {\vphantom {\left\| {\pmb \varphi
}^{k+1} -{\pmb \varphi }^{k} \right\| _{G}  \sqrt{\tilde{B}} }}
\right. \kern-\nulldelimiterspace} \sqrt{\tilde{B}} } ,
\end{equation}
we come to inequalities
\begin{equation} \label{eq_PII_58}
\sqrt{\tilde{B}} \left\| {\pmb \varphi }^{k+1} -{\pmb \varphi
}^{k} \right\| _{V_{0}} \le \left\| {\pmb \varphi }^{k+1} - {\pmb
\varphi }^{k} \right\| _{G} \le \frac{\left|\gamma \right|M_{*}
}{\sqrt{\tilde{B}} } \left\| {\pmb \varphi }^{k} \right\| _{V_{0}
}.
\end{equation}

Now let us take ${\bf v}:={\pmb \varphi }^{k+1} +{\pmb \varphi
}^{k}$ in expression (\ref {eq_PII_56}). Then $G\, ({\pmb \varphi
}^{k+1} -{\pmb \varphi }^{k} ,{\pmb \varphi }^{k+1} +{\pmb \varphi
}^{k} )=-\gamma \left[A\, ({\pmb \varphi }^{k} ,{\pmb \varphi
}^{k+1} +{\pmb \varphi }^{k} )+H_{\theta } (\bar{\bf u},{\pmb
\varphi }^{k} ,{\pmb \varphi }^{k+1} +{\pmb \varphi }^{k}
)\right]$. This relation can be written as
\[{\left\| {\pmb \varphi }^{k} \right\| _{G}^{2} -\left\|
 {\pmb \varphi }^{k+1} \right\| _{G}^{2} =\gamma \, [ 2A\, ({\pmb \varphi }^{k},
{\pmb \varphi }^{k} )+2H_{\theta } (\bar{\bf u},{\pmb \varphi
}^{k} , {\pmb \varphi }^{k} )+}\]
\[+A\, ({\pmb \varphi }^{k} ,{\pmb \varphi }^{k+1} - {\pmb \varphi
}^{k} ) +H_{\theta } (\bar{\bf u},{\pmb \varphi }^{k} ,{\pmb
\varphi }^{k+1} -{\pmb \varphi }^{k})].\]

In view of properties (\ref {eq_PII_19}) and (\ref {eq_PII_55}),
we come to an inequality
\[A\, ({\pmb \varphi }^{k} ,{\pmb \varphi }^{k+1} - {\pmb \varphi }^{k} )+
H_{\theta } (\bar{\bf u},{\pmb \varphi }^{k} ,{\pmb \varphi
}^{k+1} - {\pmb \varphi }^{k} )\ge -M_{*} \left\| {\pmb \varphi
}^{k} \right\| _{V_{0} } \left\| {\pmb \varphi }^{k+1} - {\pmb
\varphi }^{k} \right\| _{V_{0} } .\]

Suppose that $\gamma \ge 0$. Then, taking into account the
coercivity of bilinear form (\ref {eq_PII_20}), property (\ref
{eq_PII_54}), and the previous inequality, we obtain
\[\left\| {\pmb \varphi }^{k} \right\| _{G}^{2} -\left\| {\pmb \varphi }^{k+1}
 \right\| _{G}^{2} \ge \gamma \left[2B\left\| {\pmb \varphi }^{k} \right\| _{V_{0} }^{2}
 - M_{*} \left\| {\pmb \varphi }^{k} \right\| _{V_{0} } \left\| {\pmb \varphi }^{k+1} -
{\pmb \varphi }^{k} \right\| _{V_{0} } \right].\]

In view of inequalities (\ref {eq_PII_57}) and (\ref {eq_PII_58}),
we find further
\begin{equation} \label{eq_PII_59}
\left\| {\pmb \varphi }^{k} \right\| _{G}^{2} -\left\| {\pmb
\varphi }^{k+1} \right\| _{G}^{2} \ge \frac{\gamma }{\tilde{M}}
\left(2B-{\gamma M_{*}^{2} \mathord{\left/ {\vphantom {\gamma
M_{*}^{2}  \tilde{B}}} \right. \kern-\nulldelimiterspace}
\tilde{B}} \right)\, \, \left\| {\pmb \varphi }^{k} \right\|
_{G}^{2}.
\end{equation}
If the following inequality
\begin{equation} \label{eq_PII_60}
\gamma \left(2B-\gamma {M_{*}^{2} \mathord{\left/ {\vphantom
{M_{*}^{2}  \tilde{B}}} \right. \kern-\nulldelimiterspace}
\tilde{B}} \right)\, >0
\end{equation}
holds, then the sequence $\left\| {\pmb \varphi }^{k} \right\|
_{G}^{2} $ will be monotonically nonincreasing: $\left\| {\pmb
\varphi }^{k} \right\| _{G}^{2} \ge \left\| {\pmb \varphi }^{k+1}
\right\| _{G}^{2} $, and, hence, $\left\| {\pmb \varphi }^{k}
\right\| _{G}^{2} \mathop{\to }\limits_{k\to \infty } \omega $,
where $\omega \ge 0$. Passing to the limit as $k\to \infty $ in
expression (\ref {eq_PII_59}), we obtain $0\ge \frac{\gamma
}{\tilde{M}} \left(2B-{\gamma M_{*}^{2} \mathord{\left/ {\vphantom
{\gamma M_{*}^{2}  \tilde{B}}} \right. \kern-\nulldelimiterspace}
\tilde{B}} \right)\omega $, i.e. $\omega =0$, and, therefore
$\left\| {\pmb \varphi }^{k} \right\| _{G} \mathop{\to
}\limits_{k\to \infty } 0$. From inequality (\ref {eq_PII_60}) we
establish the interval of allowable values of the iterative
parameter $\gamma $:
\[\gamma \in (0\, ;\, \, \gamma _{2} ), \,\,\, \gamma _{2} ={2B\tilde{B}\mathord{\left/
{\vphantom {2B\tilde{B} M_{*}^{2} }} \right.
\kern-\nulldelimiterspace} M_{*}^{2} } .\]

Since the norms $\left\| \, \cdot \, \right\| _{G} $ and $\left\|
\, \cdot \, \right\| _{V_{0} } $ are equivalent, we have $\left\|
{\pmb \varphi }^{k} \right\| _{V_{0} } \mathop{\to }\limits_{k\to
\infty } 0$, and, hence, $\left\| {\bf u}^{k} -\bar{\bf u}\right\|
_{V_{0} } \mathop{\to }\limits_{k\to \infty } 0$.

Using inequality (\ref {eq_PII_59}), we find an estimate
\begin{equation} \label{eq_PII_61}
\left\| {\pmb \varphi }^{k+1} \right\| _{G}^{2} \le q^{2} \left\|
{\pmb \varphi }^{k} \right\| _{G}^{2} , \,\,\, q^{2} (\gamma
)=1-\frac{2B}{\tilde{M}} \gamma +\frac{M_{*}^{2}
}{\tilde{M}\tilde{B}} \gamma ^{2} .
\end{equation}

It is not hard to show, that $q^{2} \in (0;1)$ for $\gamma \in
(0\, ;\, \, \gamma _{2} )$. We obtain this from the next
relations:
\[q^{2} (0)=q^{2} (\gamma _{2} )=1, \,\,\,
 \bar{\gamma }=\arg \mathop{\min }\limits_{\gamma \in (0;\gamma _{2} )} q^{2}
(\gamma )={\gamma _{2} \mathord{\left/ {\vphantom {\gamma _{2}
2}} \right.
 \kern-\nulldelimiterspace} 2} ={B\tilde{B}\mathord{\left/
 {\vphantom {B\tilde{B} M_{*}^{2} }} \right. \kern-\nulldelimiterspace} M_{*}^{2} } .\]

As follows from estimate (\ref{eq_PII_61}), the convergence rate
is maximal if the parameter $q$ is minimal, i.e., if $\gamma
=\bar{\gamma }$. $\Box$

\textbf{Remark~1.} \textit{Suppose that the term $J'_{\theta }
({\bf u}, {\bf v})$ is G\^{a}teaux differentiable in $\bf u$. Then
conditions (\ref{eq_PII_29}), (\ref{eq_PII_30}) in
}\textbf{\textit{Theorem~4}}\textit{ can be replaced by the
following properties}
\begin{equation} \label{eq_PII_62}
\left(\forall {\bf u}, {\bf v}\in V_{0} \right)\left\{\,
J''_{\theta } ({\bf u},{\bf v},{\bf v})\ge 0\right\},
\end{equation}
\begin{equation} \label{eq_PII_63}
\left(\exists D>0\right)\left(\forall {\bf u}, {\bf v},{\bf w} \in
V_{0} \right)\left\{\, \left|J''_{\theta } ({\bf u},{\bf v},{\bf
w})\right|\le D\left\|{\bf v}\right\| _{V_{0} } \left\| {\bf
w}\right\| _{V_{0} } \right\}.
\end{equation}

\textbf{Proof.} Let us apply to $J''_{\theta} ({\bf u},{\bf
v},{\bf w})$ the Lagrange formula of finite increments \cite
{PII_4}:
\[\left(\forall {\bf u}, {\bf v},{\bf w} \in V_{0}
\right)\left(\exists \tau \in (0;1)\right)\left\{\, J''_{\theta }
({\bf u}+\tau {\bf w},{\bf w},{\bf v})=J'_{\theta } ({\bf u}+{\bf
w},{\bf v})-J'_{\theta } ({\bf u}, {\bf v})\right\}.\] Then the
satisfaction of conditions (\ref{eq_PII_62}) and (\ref{eq_PII_63})
yields the satisfaction of properties (\ref {eq_PII_30}) and (\ref
{eq_PII_31}). $\Box$

Now let us investigate the stability of iterative method (\ref
{eq_PII_46}) to the computational errors. Let us show that this
method has stability properties, which are natural for all
successive iteration methods.

Suppose that the conditions of \textbf{Theorem~4} are satisfied.
Then for any ${\bf u}^{k} \in V_{0} $, $k=0,1,...$ there exists a
unique solution ${\bf u}={\bf u}^{k+1} \in V_{0} $ of problem
(\ref {eq_PII_51}). Therefore, there exists an operator
$\mathfrak{R}:\,\, {\bf u}^{k} \in V_{0} \to {\bf u}^{k+1} \in
V_{0} $, which maps every element ${\bf u}^{k} \in V_{0} $ onto
the solution ${\bf u}={\bf u}^{k+1} $ of problem (\ref
{eq_PII_51}), and the iterative method (\ref {eq_PII_46}) can be
written in the following form
\begin{equation} \label{eq_PII_64}
{\bf u}^{k+1} =\mathfrak{R}\,({\bf u}^{k} ), \,\,\, k=0,1,...\,.
\end{equation}

Now assume that in each step $k$ of the iterative method (\ref
{eq_PII_64}) we get some computational errors. Then this iterative
method will take the form:
\begin{equation} \label{eq_PII_65}
\breve{\bf u}^{0} ={\bf u}^{0} +{\pmb \varepsilon }^{0} ,
\end{equation}
\begin{equation} \label{eq_PII_66}
\breve{\bf u}^{k+1} =\mathfrak{R}\,(\breve{\bf u}^{k} )+{\pmb
\varepsilon }^{k+1}, \,\,\, k=0,1,...\,,
\end{equation}
where $\breve{\bf u}^{k+1} $, $k=0,1,...$ is an approximate
solution of problem (\ref {eq_PII_51}), ${\pmb \varepsilon }^{k+1}
$, $k=0,1,...$ is a computational error, which occurs in each step
$k$, and ${\pmb \varepsilon }^{0} $ is an error of the initial
approximation.

\textbf{Corollary~1.} \textit{Suppose that the conditions of
}\textbf{\textit{Theorem~4}} \textit{are satisfied, and the errors
which occur in each step $k$ of the iterative method (\ref
{eq_PII_46}) are uniformly bounded, i.e.}

\[\left(\exists \varepsilon >0\right)\, \, \left(\forall k\in \left\{0,1,...\right\}\right)\, \, \left\{\, \left\| {\pmb \varepsilon }^{k} \right\| _{G} \le {\varepsilon} \, \right\}.\]
\textit{Then the following estimates hold:}

\begin{equation} \label{eq_PII_68}
\left\| \breve{\bf u}^{k} -\bar{\bf u}\right\| _{G} \le q^{k}
\left\| {\bf u}^{0} -\bar{\bf u}\right\| _{G} +\frac{\varepsilon
}{1-q}\,,
\end{equation}

\begin{equation} \label{eq_PII_69}
\left\| \breve{\bf u}^{k} -\bar{\bf u}\right\| _{G} \le
\frac{q}{1-q} \left\| \breve{\bf u}^{k} -\breve{\bf u}^{k-1}
\right\| _{G} +\frac{\varepsilon }{1-q}\,,
\end{equation}
\textit{where $\bar{\bf u}\in V_{0}$ is the exact solution of
problem (\ref {eq_PII_41}).}

The proof of this proposition follows from property (\ref
{eq_PII_50}).

Thus, from inequality (\ref {eq_PII_68}) it follows, that the
errors, which occur in each step of the iterative method (\ref
{eq_PII_46}) do not accumulate, and the number of iterations
depends linearly on the logarithm of accuracy of the initial
approximation.

Now consider a nonstationary iterative method for solution of the
nonlinear variational equation (\ref {eq_PII_41}), where bilinear
forms $G\,({\bf u},{\bf v})$ are different in each iteration \cite
{PII_47,PII_25}.

In space $V_{0}$ introduce a sequence of bilinear forms
$\left\{G^{k} :\, \, \, V_{0} \times V_{0} \to {\mathbb
R}\right\}$, $k=0,1,...$, which satisfy the property
\[\left(\forall \, Y\in V_{0}^{*} \right)\left(\exists \, !\, \,
 \bar{\bf u}\in V_{0} \right)\left(\forall \, {\bf v}\in V_{0} \right)\left\{G^{k}
 (\bar{\bf u},{\bf v})-Y({\bf v})\equiv 0\right\}.\]

For the solution of the nonlinear variational equation (\ref
{eq_PII_41}), we proposed the following nonstationary iterative
method \cite {PII_47,PII_25}:
\begin{equation} \label{eq_PII_71}
G^{k} ({\bf u}^{k+1} ,{\bf v})=G^{k} ({\bf u}^{k} ,{\bf v})-\gamma
\left[A\, ({\bf u}^{k} ,{\bf v})+J'_{\theta } ({\bf u}^{k} ,{\bf
v})-L\, (\bf{v})\right],\,\, k=0,1,...\,,
\end{equation}
where ${\bf u}^{k} \in V_{0} $ is the $k$-th approximation to the
exact solution of problem (\ref {eq_PII_41}), and $\gamma \in
{\mathbb R}$ is an iterative parameter.

We have proved the next proposition on the convergence of this
method.

\textbf{Theorem~5.} \textit{Suppose that conditions (\ref
{eq_PII_19})~--~(\ref {eq_PII_21}) and (\ref {eq_PII_29})~--~(\ref
{eq_PII_31}) hold, the bilinear forms $G^{k} ({\bf u}, {\bf v})$
satisfy the following properties}

\begin{equation} \label{eq_PII_72}
\left(\forall {\bf u}, {\bf v}\in V_{0} \right)\left\{G^{k} ({\bf
u}, {\bf v})=G^{k} ({\bf v}, {\bf u})\right\},
\end{equation}
\begin{equation} \label{eq_PII_73}
\left(\exists \tilde{M}>0\right)\left(\forall k \in
\left\{0,1,...\right\}
 \right)\left(\forall {\bf u}, {\bf v}\in
V_{0} \right)\left\{\, \left|G^{k} ({\bf u}, {\bf v})\right|\le
\tilde{M}\left\|{\bf u}\right\| _{V_{0} } \left\|{\bf v}\right\|
_{V_{0} } \right\},
\end{equation}
\begin{equation} \label{eq_PII_74}
\left(\exists \tilde{B}>0\right)\left(\forall k\in
\left\{0,1,...\right\}\right)\left(\forall {\bf u} \in V_{0}
\right)\left\{G^{k} ({\bf u}, {\bf u})\ge \tilde{B}\left\|{\bf
u}\right\| _{V_{0} }^{2} \right\},
\end{equation}
\begin{equation} \label{eq_PII_75}
\left(\exists k_{0} \in \left\{0,1,...\right\}\right)\left(\forall
k\ge k_{0} \right)\left(\forall {\bf u} \in V_{0}
\right)\left\{G^{k} ({\bf u}, {\bf u})\ge G^{k+1} ({\bf u}, {\bf
u})\right\},
\end{equation}
\textit{and the iterative parameter $\gamma$ lies in the interval
$\gamma \in (0;\gamma _{2} )$}, $\gamma _{2}
={2B\tilde{B}\mathord{\left/ {\vphantom {2B\tilde{B} (M+D)^{2} }}
\right. \kern-\nulldelimiterspace} (M+D)^{2} }$.\textit{}
\textit{Then the sequence $\{ {\bf u}^{k} \},$ obtained by the
nonstationary iterative method (\ref {eq_PII_71}), converges
strongly in $V_{0}$ to the exact solution $\bar{\bf u}\in V_{0}$
of the variational equation (\ref {eq_PII_41}), i.e. $\left\| {\bf
u}^{k} -\bar{\bf u}\right\| _{V_{0} } \mathop{\to }\limits_{k\to
\infty } 0$}.\textit{}

The proof of this theorem is similar to the proof of
\textbf{Theorem~4}. We omit the details.

Finally, let us say that it is possible to make a modification of
iterative methods (\ref {eq_PII_46}) and (\ref {eq_PII_71}) in
which the iterative parameter $\gamma$ is taken differently in
each iteration. The purpose of the nonstationary choice of
$\gamma$ might be the improvement of the convergence rate of
iterative methods (\ref {eq_PII_46}) and (\ref {eq_PII_71}).

\noindent

\section{Parallel domain decomposition schemes}

Note, that in the most general case the iterative methods (\ref
{eq_PII_46}) and (\ref {eq_PII_71}) applied to solve the nonlinear
penalty variational equations (\ref {eq_PII_41}) for multibody
contact problems do not lead to the domain decomposition.
Therefore, we now consider such variants of these methods, which
lead to the domain decomposition, namely, which reduce the
solution of original multibody contact problem in $\Omega$ to the
solution of a sequence of separate linear variational problems in
the subdomains $\Omega_{\alpha}$, $\alpha =1,2,...,N$.

Let us take the bilinear form $G$ in the iterative method (\ref
{eq_PII_46}) as follows \cite {PII_8,PII_47}:
\begin{equation} \label{eq_PII_76}
G\, ({\bf u}, {\bf v})=A\, ({\bf u}, {\bf v})+X\, ({\bf u}, {\bf
v}), \,\,\, {\bf u},{\bf v}\in V_{0} ,
\end{equation}
where $X\, ({\bf u}, {\bf v}):\, \,  V_{0} \times V_{0} \to
{\mathbb R}$ is the next bilinear form \cite {PII_8,PII_47}:
\begin{equation} \label{eq_PII_77}
X({\bf u}, {\bf v})=\frac{1}{\theta} \sum _{\left\{\alpha ,\, \,
\beta \right\}\, \in \, Q}\int _{\, S_{\alpha \beta }
}\left(u_{\alpha \, n} v_{\alpha \, n} \, \psi _{\alpha \beta }
+u_{\beta \, n} v_{\beta \, n} \, \psi _{\beta \alpha } \right)\,
dS , \,\,\, {\bf u},{\bf v}\in V_{0} .
\end{equation}
Here $\psi _{\alpha \beta } ({\bf x})=\left\{0,\, \, \, {\bf x}
\in S_{\alpha \beta } \backslash S_{\alpha \beta }^{1}
\right\}\vee \left\{1,\, \, \, {\bf x} \in S_{\alpha \beta }^{1}
\right\}$ are the characteristic functions of some given subareas
$S_{\alpha \beta }^{1} \subseteq S_{\alpha \beta } $ of the
possible contact zones $S_{\alpha \beta } $, $\alpha =1,2,...,N$,
$\beta \in B_{\alpha } $.

Introduce a notation $\tilde{\bf u}^{k+1}=\frac{1}{\gamma}[{\bf
u}^{k+1}-(1-\gamma){\bf u}^{k}]$. Then the iterative method (\ref
{eq_PII_46}) with bilinear form (\ref {eq_PII_77}) can be written
in the following equivalent way:
\begin{equation} \label{eq_PII_78}
A \left(\tilde{\bf u}^{k+1} ,{\bf v}\right)+X \left(\tilde{\bf
u}^{k+1} ,{\bf v}\right)=L\, ({\bf v})+X({\bf u}^{k} ,{\bf
v})-J'_{\theta } ({\bf u}^{k} ,{\bf v}), \,\,\, \forall \, {\bf
v}\in V_{0} ,
\end{equation}
\begin{equation} \label{eq_PII_79}
{\bf u}^{k+1} =\gamma \, \tilde{\bf u}^{k+1} +\left(1-\gamma
\right){\bf u}^{k} , \,\,\, k=0,1,...\,.
\end{equation}

\textbf{Lemma~4.} \textit{Suppose that the surfaces $S_{\alpha
\beta}$},\textit{ $\left\{ \alpha, \beta \right\} \in Q,$ are
Lipschitz. Then the bilinear form (\ref {eq_PII_77}) is symmetric,
continuous and nonnegative, i.e.}

\begin{equation} \label{eq_PII_80}
\left(\forall {\bf u}, {\bf v}\in V_{0} \right)\left\{X\, ({\bf
u}, {\bf v})=X\, ({\bf v}, {\bf u})\right\},
\end{equation}

\begin{equation} \label{eq_PII_81}
\left(\exists Z>0\right)\left(\forall {\bf u}, {\bf v}\in V_{0}
\right)\left\{\, \left|X\, ({\bf u}, {\bf v})\right|\le
Z\left\|{\bf u}\right\| _{V_{0} } \left\|{\bf v}\right\| _{V_{0} }
\right\},
\end{equation}

\begin{equation} \label{eq_PII_82}
\left(\forall {\bf u} \in V_{0} \right)\left\{X\, ({\bf u}, {\bf
u})\ge 0\right\}.
\end{equation}

\textbf{Proof.} It is obvious that conditions (\ref {eq_PII_80})
and (\ref {eq_PII_82}) hold. Thus, let us show the continuity of
bilinear form (\ref {eq_PII_77}).

We can write $X\, ({\bf u}, {\bf v})=\sum _{\left\{\alpha ,\, \,
\beta \right\}\, \in \, Q}X_{\alpha \beta } ({\bf u}, {\bf v})$,
where
\[ X_{\alpha \beta } ({\bf u}, {\bf v})=\frac{1}{\theta }
\left(\int _{\, S_{\alpha \beta } }\psi _{\alpha \beta } \,
u_{\alpha \, n} v_{\alpha \, n} \, dS +\int _{\, S_{\alpha \beta }
}\psi _{\beta \alpha } \, u_{\beta \, n} v_{\beta \, n} \, dS
\right), \,\,\, \left\{\alpha ,\beta \right\}\in Q. \]

The first term can be written in the following way:
\[\int _{\, S_{\alpha \beta } }\psi _{\alpha \beta } \, u_{\alpha
\, n} v_{\alpha \, n} \, dS =\int _{S_{\alpha \beta } }\psi
_{\alpha \beta } \, \left[{\bf n}_{\alpha } \cdot
{\rm{Tr}}_{\alpha }^{0} ({\bf u}_{\alpha}) \right] \left[{\bf
n}_{\alpha } \cdot {\rm{Tr}}_{\alpha }^{0} ({\bf v}_{\alpha})
\right]\, dS.  \] Taking into account that the functions $\psi
_{\alpha \beta } $ and the components of unit normals ${\bf
n}_{\alpha } $ are bounded, and using inequality (\ref
{eq_PII_33}) and the Schwarz inequality, we obtain
\[\left(\int _{S_{\alpha \beta } }\psi _{\alpha \beta } \,
 u_{\alpha \, n} v_{\alpha \, n} \, dS \right)^{2} \le 9\left\| {\rm{Tr}}_{\, \alpha }^{0} ({\bf u}_{\alpha } )\right\| _{[L_{2}
 (S_{\alpha \beta } )]^{3} }^{2} \left\| {\rm{Tr}}_{\, \alpha }^{0}
 ({\bf v}_{\alpha } )\right\| _{[L_{2} (S_{\alpha \beta } )]^{3} }^{2} .\]
In view of inequality (\ref {eq_PII_35}), we find further

\[\left|\int _{S_{\alpha \beta } }\psi _{\alpha \beta } \,
u_{\alpha \, n} v_{\alpha \, n} \, dS \right|\le 3\,
\tilde{W}_{\alpha \beta }
 \left\| {\bf u}_{\alpha } \right\| _{V_{\alpha } } \left\| {\bf v}_{\alpha }
  \right\| _{V_{\alpha } } ,\,\, \tilde{W}_{\alpha \beta } =T_{\alpha
  \beta }^{2} +\tilde{\varepsilon }_{\alpha \beta } >0, \,
  \,\tilde{\varepsilon }_{\alpha \beta } >0.\]

Similar inequality can be obtained for the second term of
$X_{\alpha \beta}$. Thus,
\[\left|X_{\alpha \beta } ({\bf u}, {\bf v})\right|\le \frac{3\,
 W_{\alpha \beta } }{\theta } \left(\left\| {\bf u}_{\alpha }
  \right\| _{V_{\alpha } } \left\| {\bf v}_{\alpha } \right\| _{V_{\alpha } }
   +\left\| {\bf u}_{\beta } \right\| _{V_{\beta } } \left\| {\bf v}_{\beta }
   \right\| _{V_{\beta } } \right), \,\,
    W_{\alpha \beta } =\max \left\{\tilde{W}_{\alpha \beta },\,
    \tilde{W}_{\beta \alpha } \right\}.\]
As a result, we establish
\[\left|X({\bf u}, {\bf v})\right|\le \frac{3\, W}{\theta }
\sum _{\alpha =1}^{N}\sum _{\beta =1}^{N}\left(\left\| {\bf
u}_{\alpha } \right\| _{V_{\alpha } } \left\| {\bf v}_{\alpha }
\right\| _{V_{\alpha } } +\left\| {\bf u}_{\beta } \right\|
_{V_{\beta } } \left\| {\bf v}_{\beta } \right\| _{V_{\beta } }
\right)  =\]

\[=\frac{6\, WN}{\theta } \sum _{\alpha =1}^{N}\left\| {\bf u}_{\alpha } \right\| _{V_{\alpha } } \left\| {\bf v}_{\alpha } \right\| _{V_{\alpha } }
\le Z\left\|{\bf u}\right\| _{V_{0} } \left\|{\bf v}\right\|
_{V_{0} } , \,\,\, \forall \, {\bf u}, {\bf v}\in V_{0} ,\] where
$Z={6\, WN\mathord{\left/ {\vphantom {6\, WN \theta }} \right.
\kern-\nulldelimiterspace} \theta } >0$, $W=\mathop{\max
}\limits_{1\le \alpha ,\beta \le N} W_{\alpha \beta }>0$. $\Box$

From \textbf{Lemma~4} and \textbf{Lemma~1}, it follows that the
bilinear form (\ref {eq_PII_76}) is symmetric, continuous and
coercive with constants $\tilde{M}=M+Z$ and $\tilde{B}=B$
respectively. In addition, due to \textbf{Lemmas} \textbf{1},
\textbf{2}, and \textbf{4}, the functionals $L\,(\bf{v})$, $X({\bf
u}^{k},{\bf v})$, and $J'_{\theta} ({\bf u}^{k},{\bf v})$ are
linear and continuous in ${\bf v}$. Therefore, there exists a
unique solution ${\bf u}=\tilde{\bf u}^{k+1} \in V_{0} $ of the
variational problem (\ref {eq_PII_78}).

Thus, the conditions of \textbf{Theorem~4} are satisfied, and we
obtain the next proposition.

\textbf{Theorem~6.} \textit{Suppose that the conditions of
\textbf{Lemmas~1} and \textbf{2} hold, and  $\gamma \in (0;\gamma
_{2})$}, $\gamma_{2} ={2B^{2} \mathord{\left/ {\vphantom {2B^{2}
(M+D)^{2}}} \right. \kern-\nulldelimiterspace} (M+D)^{2}}$.
\textit{Then the sequence $\{ {\bf u}^{k} \},$} \textit{obtained
by the iterative method (\ref {eq_PII_78})~--~(\ref {eq_PII_79}),
which is equivalent to the iterative method (\ref {eq_PII_46})
with bilinear form (\ref {eq_PII_76}), converges strongly in
$V_{0} $ to the exact solution $\bar{\bf u}\in V_{0} $ of the
nonlinear penalty variational equation (\ref {eq_PII_41}) for the
unilateral multibody contact problem, i.e.} $\left\| {\bf u}^{k}
-\bar{\bf u}\right\| _{V_{0} } \mathop{\to }\limits_{k\to \infty }
0$. \textit{Moreover, the convergence rate in the norm $\left\| \,
\cdot \, \right\| _{G} $ is linear (\ref {eq_PII_50}), where
$q=\sqrt{1-{\gamma \left(2B-\gamma M_{*}^{2}
/B\right)\mathord{\left/ {\vphantom {\gamma \left(2B-\gamma
M_{*}^{2} /B\right) \left(M+Z\right)}} \right.
\kern-\nulldelimiterspace} \left(M+Z\right)} \, } $}, \textit{and
the maximal rate reaches as $\gamma =\bar{\gamma }={B^{2}
\mathord{\left/ {\vphantom {B^{2}  M_{*}^{2} }} \right.
\kern-\nulldelimiterspace} M_{*}^{2} } $}, $M_{*} =M+D$.

Now let us show that the iterative method (\ref
{eq_PII_78})~--~(\ref {eq_PII_79}) leads to the domain
decomposition.

Due to the relationships
\[J'_{\theta } ({\bf u}, {\bf v})=-\frac{1}{\theta } \sum _{\alpha =1}^{N}\sum _{\,
 \beta \, \in B_{\alpha } }\int _{\, S_{\alpha \beta } }
 \left(d_{\alpha \beta } -u_{\alpha \, n} -u_{\beta \, n} \right)^{-}
 v_{\alpha \, n} \, dS,\]
\[X({\bf u}, {\bf v})=\frac{1}{\theta } \sum _{\alpha
=1}^{N}\sum _{\, \beta \, \in B_{\alpha } }\int _{\, S_{\alpha
\beta } }\psi _{\alpha \beta } \, u_{\alpha \, n} \, v_{\alpha \,
n} \, dS,\] method (\ref {eq_PII_78})~--~(\ref {eq_PII_79})
rewrites as follows:

\[\sum _{\alpha =1}^{N}a_{\alpha } (\tilde{\bf u}_{\alpha }^{k+1} ,{\bf v}_{\alpha } )
+\frac{1}{\theta } \sum _{\alpha =1}^{N}\sum _{\, \beta \, \in
B_{\alpha } }\int _{\, S_{\alpha \beta } }\psi _{\alpha \beta }
 \left(\tilde{u}_{\alpha \, n}^{k+1} -u_{\alpha \, n}^{k} \right)v_{\alpha \, n} \, dS   =\]

\begin{equation} \label{eq_PII_83}
=\sum _{\alpha =1}^{N}l_{\alpha } ({\bf v}_{\alpha } )
+\frac{1}{\theta } \sum _{\alpha =1}^{N}\sum _{\, \beta \, \in
B_{\alpha } }\int _{\, S_{\alpha \beta } }\left(d_{\alpha \beta }
-u_{\alpha \, n}^{k} -u_{\beta \, n}^{k} \right)^{-} v_{\alpha \,
n} \, dS   ,
\end{equation}

\begin{equation} \label{eq_PII_84}
{\bf u}_{\alpha }^{k+1} =\gamma \, \tilde{\bf u}_{\alpha }^{k+1}
+\left(1-\gamma \right){\bf u}_{\alpha }^{k} , \,\,\, \alpha
=1,2,...,N, \,\,\, k=0,1,...\,.
\end{equation}

Since the common quantities of the subdomains are known from the
previous iteration, the variational equation (\ref {eq_PII_83})
splits into $N$ variational equations in the separate subdomains
$\Omega_{\alpha}$. Therefore, method (\ref {eq_PII_83})~--~(\ref
{eq_PII_84}) can be written in the following equivalent form:

\[{a_{\alpha } (\tilde{\bf u}_{\alpha }^{k+1} ,\,
\, {\bf v}_{\alpha } )+\frac{1}{\theta } \sum _{\beta \,
 \in B_{\alpha } }\int _{\, S_{\alpha \beta } }\psi _{\alpha \beta }
 \, \tilde{u}_{\alpha \, n}^{k+1} \, v_{\alpha \, n} \, dS
=}\]
\[{=l_{\alpha } ({\bf v}_{\alpha } )+\frac{1}{\theta } \sum _{\beta
\, \in B_{\alpha } }\int _{\, S_{\alpha \beta } }\psi _{\alpha
\beta } \, u_{\alpha \, n}^{k} \, v_{\alpha \, n} \, dS  +}\]

\begin{equation} \label{eq_PII_85}
+\frac{1}{\theta } \sum _{\beta \, \in B_{\alpha } }\int _{\,
S_{\alpha \beta } }\left(d_{\alpha \beta } -u_{\alpha \, n}^{k}
-u_{\beta \, n}^{k} \right)^{-} v_{\alpha \, n} \, dS, \,\,\,
\forall \, {\bf v}_{\alpha } \in V_{\alpha }^{0} ,
\end{equation}

\begin{equation} \label{eq_PII_86}
{\bf u}_{\alpha }^{k+1} =\gamma \, \tilde{\bf u}_{\alpha }^{k+1}
+\left(1-\gamma \right){\bf u}_{\alpha }^{k} , \,\,\, \alpha
=1,2,...,N, \,\,\, k=0,1,...\,.
\end{equation}

Since the bilinear forms $a_{\alpha } ({\bf u}_{\alpha } ,{\bf
v}_{\alpha } )$, $X_{\alpha } ({\bf u}_{\alpha } ,{\bf v}_{\alpha
} )=\frac{1}{\theta } \sum _{\, \beta \, \in B_{\alpha } }\int
_{\, S_{\alpha \beta } }\psi _{\alpha \beta } \, u_{\alpha \, n}
v_{\alpha \, n} \, dS $ are symmetric, continuous and coercive,
and the functionals $l_{\alpha } ({\bf v}_{\alpha } )$, $X_{\alpha
} ({\bf u}_{\alpha }^{k} ,{\bf v}_{\alpha } )$, $\frac{1}{\theta }
\sum _{\beta \, \in B_{\alpha } }\int _{\, S_{\alpha \beta }
}\left(d_{\alpha \beta } -u_{\alpha \, n}^{k} -u_{\beta \, n}^{k}
\right)^{-} v_{\alpha \, n} \, dS  $ are linear and continuous in
${\bf v}_{\alpha } $, it follows that there exists a unique
solution ${\bf u}_{\alpha }^{*} =\tilde{\bf u}_{\alpha }^{k+1} \in
V_{\alpha }^{0}$ of each variational equation (\ref {eq_PII_85}).
Furthermore, it is obvious to see that a unique solution ${\bf
u}^{*}$ of the variational equation (\ref {eq_PII_78}), i.e. (\ref
{eq_PII_83}), takes the form ${\bf u}^{*} =\tilde{\bf u}^{k+1}
=\left({\bf u}_{1}^{*},{\bf u}_{2}^{*},...,{\bf u}_{N}^{*}
\right)^{{\rm T}} \in V_{0} $. Therefore, the solution of the
variational equation (\ref {eq_PII_78}) is equivalent to the
solution of $N$ variational equations (\ref {eq_PII_85}) in the
separate subdomains, and the iterative processes (\ref
{eq_PII_78})~--~(\ref {eq_PII_79}) and (\ref {eq_PII_85})~--~(\ref
{eq_PII_86}) are equivalent.

Now, consider the iterative method (\ref {eq_PII_85})~--~(\ref
{eq_PII_86}) in more detail.

In each iteration of this method we have to solve $N$ linear
variational equations (\ref {eq_PII_85}) in parallel, which
correspond to some linear elasticity problems in subdomains with
prescribed Robin boundary conditions on the possible contact
areas:
\begin{equation} \label{eq_PII_87}
\tilde{\sigma }_{\alpha \beta n}^{k+1} +\psi _{\alpha \beta }
{\tilde{u}_{\alpha \, n}^{k+1} \mathord{\left/ {\vphantom
{\tilde{u}_{\alpha \, n}^{k+1}  \theta }} \right.
\kern-\nulldelimiterspace} \theta } ={\left(d_{\alpha \beta }
-u_{\alpha \, n}^{k} -u_{\beta \, n}^{k} \right)^{-}
\mathord{\left/ {\vphantom {\left(d_{\alpha \beta } -u_{\alpha \,
n}^{k} -u_{\beta \, n}^{k} \right)^{-}  \theta }} \right.
\kern-\nulldelimiterspace} \theta } +{\psi _{\alpha \beta } \,
u_{\alpha \, n}^{k} \mathord{\left/ {\vphantom {\psi _{\alpha
\beta } \, u_{\alpha \, n}^{k}  \theta }} \right.
\kern-\nulldelimiterspace} \theta } \,\,\,\, {\rm on} \,\,\,
S_{\alpha \beta}.
\end{equation}
Here $\tilde{\sigma}_{\alpha \beta n}^{k+1}$ are unknown normal
stresses on the possible contact areas $S_{\alpha \beta}$.
Therefore, the iterative method (\ref {eq_PII_85})~--~(\ref
{eq_PII_86}) refers to the parallel Robin--Robin type domain
decomposition schemes.

Since the domain decomposition method (\ref {eq_PII_85})~--~(\ref
{eq_PII_86}) and the iterative method (\ref {eq_PII_78})~--~(\ref
{eq_PII_79}) are equivalent, the convergence \textbf{Theorem~6}
also holds for method (\ref {eq_PII_85})~--~(\ref {eq_PII_86}).

Note, that from \textbf{Theorem~6}, it follows that the domain
decomposition method (\ref {eq_PII_85})~--~(\ref {eq_PII_86}) is
convergent for arbitrary initial approximations ${\bf
u}_{\alpha}^{0}$, $\alpha=1,2,...,N$, and its convergence rate is
linearly depended on them.

Taking different characteristic functions $\psi_{\alpha \beta}$ in
(\ref {eq_PII_85}), i.e. different subareas $S_{\alpha \beta}^{1}$
of the possible contact zones $S_{\alpha \beta}$, we can obtain
different particular cases of the domain decomposition method
(\ref {eq_PII_85})~--~(\ref {eq_PII_86}).

Thus, taking $\psi _{\alpha \beta } ({\bf x})\equiv 1$, $\forall
\alpha ,\beta $, i.e. $S_{\alpha \beta }^{1} =S_{\alpha \beta } $,
we get a domain decomposition scheme with Robin boundary
conditions on the whole possible contact areas:
\[\tilde{\sigma }_{\alpha \beta n}^{k+1} +{\tilde{u}_{\alpha \, n}^{k+1}
 \mathord{\left/ {\vphantom {\tilde{u}_{\alpha \,
 n}^{k+1}  \theta }} \right. \kern-\nulldelimiterspace} \theta }
  ={\left(d_{\alpha \beta } -u_{\alpha \, n}^{k} -u_{\beta \,
   n}^{k} \right)^{-} \mathord{\left/ {\vphantom {\left(d_{\alpha \beta }
    -u_{\alpha \, n}^{k} -u_{\beta \, n}^{k} \right)^{-}  \theta }} \right.
     \kern-\nulldelimiterspace} \theta } +{u_{\alpha \, n}^{k} \mathord{\left/
      {\vphantom {u_{\alpha \, n}^{k}  \theta }} \right. \kern-\nulldelimiterspace}
       \theta }  \,\,\,\,  {\rm on} \,\,\, S_{\alpha \beta}.\]
Therefore, we have called this domain decomposition method as full
parallel Robin--Robin domain decomposition scheme \cite {PII_8}.

Taking $\psi _{\alpha \beta } ({\bf x})\equiv 0$, $\forall \alpha
,\beta $, i.e. $S_{\alpha \beta }^{1} =\emptyset $, we get a
parallel Neumann--Neumann domain decomposition scheme \cite
{PII_23,PII_7,PII_8}:
\begin{equation} \label{eq_PII_88}
a_{\alpha } (\tilde{\bf u}_{\alpha }^{k+1} ,\, \, {\bf v}_{\alpha
} )=l_{\alpha } ({\bf v}_{\alpha } )+\frac{1}{\theta } \sum
_{\beta \, \in B_{\alpha } }\int _{\, S_{\alpha \beta }
}\left(d_{\alpha \beta } -u_{\alpha \, n}^{k} -u_{\beta \, n}^{k}
\right)^{-} v_{\alpha \, n} \, dS  ,\,\,\, \forall \, {\bf
v}_{\alpha } \in V_{\alpha }^{0} ,
\end{equation}

\begin{equation} \label{eq_PII_89}
{\bf u}_{\alpha }^{k+1} =\gamma \, \tilde{\bf u}_{\alpha }^{k+1}
+\left(1-\gamma \right){\bf u}_{\alpha }^{k} , \,\,\, \alpha
=1,2,...,N, \,\,\, k=0,1,...\,.
\end{equation}

In each step $k$ of this scheme we have to solve $N$ variational
equations (\ref {eq_PII_88}) in parallel, which correspond to
elasticity problems in subdomains with Neumann boundary conditions
on the possible contact areas:
\[\tilde{\sigma}_{\alpha \beta n}^{k+1} =\sigma_{\alpha \beta n}^{k} ={\left(d_{\alpha
\beta} -u_{\alpha \, n}^{k} -u_{\beta \, n}^{k} \right)^{-}
\mathord{\left/ {\vphantom {\left(d_{\alpha \beta} -u_{\alpha \,
n}^{k} -u_{\beta \, n}^{k} \right)^{-}  \theta }} \right.
\kern-\nulldelimiterspace} \theta } \,\,\,\, {\rm on} \,\,\,
S_{\alpha \beta}.\]

Note, that in the most general case, we can choose functions
$\psi_{\alpha \beta}$, i.e. surfaces $S_{\alpha \beta }^{1}$,
differently for each $\alpha, \beta$.

Moreover, we can choose functions $\psi_{\alpha \beta}$
differently at each iteration $k$, i.e.
\begin{equation} \label{eq_PII_90}
\psi_{\alpha \beta} ({\bf x})=\psi_{\alpha \beta}^{k} ({\bf
x})=\left\{0,\, \, \, {\bf x} \in S_{\alpha \beta} \backslash
S_{\alpha \beta}^{k} \right\}\vee \left\{1,\, \, \, {\bf x} \in
S_{\alpha \beta}^{k} \right\},
\end{equation}
where $S_{\alpha \beta}^{k} \subseteq S_{\alpha \beta}$,
$k=0,1,...$ are some given subareas of the possible contact zones
$S_{\alpha \beta}$, $\alpha =1,2,...,N$, $\beta \in B_{\alpha}$.
As a result we obtain the following nonstationary Robin--Robin
type domain decomposition scheme

\[a_{\alpha } (\tilde{\bf u}_{\alpha }^{k+1} ,\,
\, {\bf v}_{\alpha } )+\frac{1}{\theta } \sum _{\beta \,
 \in B_{\alpha } }\int _{\, S_{\alpha \beta } }\psi _{\alpha \beta }^{k} \left(\,
 \tilde{u}_{\alpha \, n}^{k+1} -u_{\alpha \, n}^{k} \right)v_{\alpha \, n} \, dS  =\]

\begin{equation} \label{eq_PII_91}
=l_{\alpha } ({\bf v}_{\alpha } )+\frac{1}{\theta } \sum _{\beta
\, \in B_{\alpha } }\int _{\, S_{\alpha \beta } }\left(d_{\alpha
\beta } -u_{\alpha \, n}^{k} -u_{\beta \, n}^{k} \right)^{-}
v_{\alpha \, n} \, dS  , \,\,\, \forall \, {\bf v}_{\alpha } \in
V_{\alpha }^{0} ,
\end{equation}

\begin{equation} \label{eq_PII_92}
{\bf u}_{\alpha }^{k+1} =\gamma \, \tilde{\bf u}_{\alpha }^{k+1}
+\left(1-\gamma \right){\bf u}_{\alpha }^{k} , \,\,\, \alpha
=1,2,...,N, \,\,\, k=0,1,...\,.
\end{equation}

This domain decomposition scheme is equivalent to the
nonstationary iterative method (\ref {eq_PII_71}) with bilinear
forms
\begin{equation} \label{eq_PII_93}
G^{k} ({\bf u}, {\bf v})=A\, ({\bf u}, {\bf v})+X^{k} ({\bf u},
{\bf v}), \,\,\, {\bf u},{\bf v}\in V_{0} , \,\,\, k=0,1,...\,,
\end{equation}
where
\begin{equation} \label{eq_PII_94}
X^{k} ({\bf u}, {\bf v})=\frac{1}{\theta } \sum _{\left\{\alpha
,\, \, \beta \right\}\, \in \, Q}\int _{\, S_{\alpha \beta }
}\left(u_{\alpha \, n} v_{\alpha \, n} \, \psi _{\alpha \beta
}^{k} +u_{\beta \, n} v_{\beta \, n} \, \psi _{\beta \alpha }^{k}
\right)\, dS , \,\, {\bf u},{\bf v}\in V_{0} ,\,\, k=0,1,...\,.
\end{equation}

If all of the surfaces $S_{\alpha \beta}$, $\left\{ \alpha, \beta
\right\} \in Q$, are Lipschitz, then bilinear forms $X^{k} ({\bf
u}, {\bf v})$ are symmetric, nonnegative, and continuous with
constant $Z>0$. Hence, bilinear forms (\ref {eq_PII_93}) satisfy
properties (\ref {eq_PII_72})~--~(\ref {eq_PII_74}), where
$\tilde{M}=M+Z$, $\tilde{B}=B$.

It is obvious to see that condition (\ref {eq_PII_75}) in
\textbf{Theorem~5} for bilinear forms (\ref {eq_PII_93}) is
equivalent to the following condition
\[\left(\exists k_{0} \in \left\{0,1,...\right\}\right)\left(\forall k\ge k_{0}
\right)\left\{X^{k} ({\bf u}, {\bf u})\ge X^{k+1} ({\bf u}, {\bf
u})\right\},\] which by-turn is equivalent to the condition
\begin{equation} \label{eq_PII_95}
\left(\exists k_{0} \in \left\{0,1,...\right\}\right)\left(\forall
k\ge k_{0} \right)\left(\forall \alpha \right)\left(\forall \beta
\in B_{\alpha } \right)\left(\forall {\bf x} \in S_{\alpha \beta }
\right)\left\{\psi _{\alpha \beta }^{k} ({\bf x})\ge \psi _{\alpha
\beta }^{k+1} ({\bf x})\right\}.
\end{equation}

Therefore, from \textbf{Theorem~5} we obtain the next proposition
on the convergence of nonstationary domain decomposition scheme
(\ref {eq_PII_91})~--~(\ref {eq_PII_92}).

\textbf{Theorem~7.} \textit{Suppose that the conditions of
\textbf{Lemmas~1,~2} hold, $\gamma \in (0;\gamma _{2} )$},
$\gamma_{2} ={2B^{2} \mathord{\left/ {\vphantom {2B^{2} (M+D)^{2}
}} \right. \kern-\nulldelimiterspace} (M+D)^{2} } $,\textit{ and
the functions $\psi_{\alpha \beta }^{k}$ satisfy property (\ref
{eq_PII_95}). Then the sequence $\{ {\bf u}^{k} \}$},\textit{
obtained by the nonstationary domain decomposition scheme (\ref
{eq_PII_91})~--~(\ref {eq_PII_92}), converges strongly in $V_{0}$
to the exact solution $\bar{\bf u} \in V_{0}$ of the nonlinear
penalty variational equation (\ref {eq_PII_41}) for the unilateral
multibody contact problem, i.e. $\left\| {\bf u}^{k} -\bar{\bf
u}\right\| _{V_{0} } \mathop{\to }\limits_{k\to \infty } 0$}.

Now let us consider a particular case of the domain decomposition
method (\ref {eq_PII_91})~--~(\ref {eq_PII_92}). In each iteration
$k$ let us choose the functions $\psi_{\alpha \beta}^{k}$ as
follows \cite {PII_23,PII_8,PII_47,PII_12}:
\begin{equation} \label{eq_PII_96}
\psi _{\alpha \beta }^{k} ({\bf x})=\chi _{\alpha \beta }^{k}
({\bf x})=\left\{\begin{array}{c} {0,\, \, \, d_{\alpha \beta }
({\bf x})-u_{\alpha \, n}^{k} ({\bf x})-u_{\beta \, n}^{k} ({\bf
x'})\ge 0}
\\ {1,\, \, \, \, d_{\alpha \beta } ({\bf x})-u_{\alpha \, n}^{k} ({\bf x})-u_{\beta \, n}^{k}
({\bf x'})<0} \end{array}\right. ,
\end{equation}
where ${\bf x} \in S_{\alpha \beta } $, ${\bf x'}=P({\bf x})\in
S_{\beta \alpha } $. Then, taking into consideration that
$\left(d_{\alpha \beta } -u_{\alpha \, n}^{k} -u_{\beta \, n}^{k}
\right)^{-} =\left(d_{\alpha \beta } -u_{\alpha \, n}^{k}
-u_{\beta \, n}^{k} \right)\chi _{\alpha \beta }^{k} $, we obtain
the method \cite {PII_23,PII_8,PII_12}:
\begin{equation} \label{eq_PII_97}
a_{\alpha } (\tilde{\bf u}_{\alpha }^{k+1} ,\, \, {\bf v}_{\alpha
} )+\frac{1}{\theta } \sum _{\beta \, \in B_{\alpha } }\int _{\,
S_{\alpha \beta } }\chi _{\alpha \beta }^{k} \,
\left(\tilde{u}_{\alpha \, n}^{k+1} -(d_{\alpha \beta } -u_{\beta
\, n}^{k} )\right)v_{\alpha \, n} \, dS  =l_{\alpha } ({\bf
v}_{\alpha } ),
\end{equation}

\begin{equation} \label{eq_PII_98}
{\bf u}_{\alpha }^{k+1} =\gamma \, \tilde{\bf u}_{\alpha }^{k+1}
+\left(1-\gamma \right){\bf u}_{\alpha }^{k} , \,\,\, \alpha
=1,2,...,N, \,\,\, k=0,1,...\,.
\end{equation}

In each step $k$ of this method we have to solve $N$ variational
equations (\ref {eq_PII_97}) in parallel, which correspond to the
elasticity problems in subdomains with prescribed displacements
$d_{\alpha \beta} - u_{\beta \, n}^{k}$ through the penalty on
some subareas of the possible contact zones $S_{\alpha \beta}$.
Therefore, we can conventionally name this method as the
nonstationary parallel Dirichlet--Dirichlet domain decomposition
scheme.

The advantages of proposed domain decomposition schemes are their
simplicity, and the regularization of the original contact problem
because of the use of the penalty term. These domain decomposition
schemes have only one iterative loop, which deals with the domain
decomposition and the nonlinearity of unilateral contact
conditions.

Presented domain decomposition algorithms can be modified to solve
more complicated contact problems. In work \cite {PII_47} we
proposed a generalization of these algorithms to the solution of
unilateral multibody contact problems of nonlinear elasticity, and
in work \cite {PII_48} we generalized them to solve the problems
of unilateral contact between elastic bodies with nonlinear
Winkler covers. In works \cite {PII_47,PII_49} we obtained the
similar penalty domain decomposition methods for ideal multibody
contact problems and showed their connection with classical DDMs
without penalty \cite {PII_13}. The primary feature of the
algorithms presented in \cite {PII_47,PII_48} is that they deal
with all of the nonlinearities of the problem (the nonlinearity of
unilateral contact conditions, the nonlinearity of stress--strain
relationship, the nonlinearity of Winkler covers) and the domain
decomposition in one iterative loop.

Finally, let us say that the iterative methods (\ref {eq_PII_46})
and (\ref {eq_PII_71}) for the solution of nonlinear variational
equations are rather general. From these methods, besides the
parallel Robin--Robin type domain decomposition schemes (\ref
{eq_PII_85})~--~(\ref {eq_PII_86}) and (\ref {eq_PII_91})~--~(\ref
{eq_PII_92}), we can also obtain other different particular
iterative methods for the solution of the penalty variational
equation of multibody unilateral contact problems, which do not
lead to the domain decomposition.

Thus, taking the bilinear form $G$ as follows
\begin{equation} \label{eq_PII_99}
G\, ({\bf u}, {\bf v})=A\, ({\bf u}, {\bf v})+\tilde{X}\, ({\bf
u}, {\bf v}),
\end{equation}
\begin{equation} \label{eq_PII_100}
\tilde{X}\, ({\bf u}, {\bf v})=\frac{1}{\theta } \sum
_{\left\{\alpha ,\, \, \beta \right\}\, \in \, Q}\int _{\,
S_{\alpha \beta } }\left(u_{\alpha \, n} +u_{\beta \, n}
\right)\left(v_{\alpha \, n} +v_{\beta \, n} \right)\, dS,
\end{equation}
and the iterative parameter $\gamma=1$ in (\ref {eq_PII_46}), we
obtain the iterative method for the solution of multibody
unilateral contact problems, which can be viewed as a
generalization of the penalty iteration method, proposed in \cite
{PII_14} for the solution of crack problems with nonpenetration
condition.

Taking the bilinear forms $G^{k}$ in each step $k$ of method (\ref
{eq_PII_71}) as follows
\begin{equation} \label{eq_PII_101}
G^{k}({\bf u},{\bf v})={\partial}^{2} F_{\theta}({\bf u}^{k},{\bf
u},{\bf v})=A\, ({\bf u}, {\bf v})+\tilde{X}^{k} ({\bf u}, {\bf
v}), \,\,\, {\bf u},{\bf v}\in V_{0},
\end{equation}
and the iterative parameter $\gamma=1$, where ${\partial}^{2}
F_{\theta}({\bf u}^{k},{\bf u},{\bf v})$ is one of the second
G\^{a}teaux subdifferentials of the functional $F_{\theta}$ in the
point ${\bf u}^{k} \in V_{0}$, and
\[
\tilde{X}^{k} ({\bf u},{\bf v})=\frac{1}{\theta} \sum
_{\left\{\alpha,\, \beta \right\}\, \in \, Q}\int _{S_{\alpha
\beta}}\left(u_{\alpha \, n} +u_{\beta \, n}
\right)\left(v_{\alpha \, n} +v_{\beta \, n} \right)\, \chi
_{\alpha \beta }^{k} \, dS ,\,\, {\bf u},{\bf v}\in V_{0},
\]
\[\chi _{\alpha \beta }^{k} ({\bf x})=\left\{\begin{array}{c} {0,\, \,
\, d_{\alpha \beta } ({\bf x})-u_{\alpha \, n}^{k}
 ({\bf x})-u_{\beta \, n}^{k} ({\bf x'})\ge 0} \\ {1,\, \, \,
 \, d_{\alpha \beta } ({\bf x})-u_{\alpha \, n}^{k} ({\bf x})-u_{\beta \,
  n}^{k} ({\bf x'})<0} \end{array}\right.,\,\,{\bf x} \in S_{\alpha \beta } ,\,
  \, {\bf x'}=P({\bf x})\in S_{\beta \alpha } ,\]
we obtain the iterative method, which can be viewed as an active
set method, i.e. an implicit semi-smooth Newton method for
unilateral multibody contact problems. The convergence theorem for
the active set method for crack problems with nonpenetration
condition, as a variant of the semi-smooth Newton method, was
proved in \cite {PII_10}.

Note, that the bilinear forms $\tilde{X}\, ({\bf u}, {\bf v})$ and
$\tilde{X}^{k} ({\bf u}, {\bf v})$ are symmetric, nonnegative and
coercive. However, the iterative methods (\ref {eq_PII_46}), (\ref
{eq_PII_99}) and (\ref {eq_PII_71}), (\ref {eq_PII_101}) do not
lead to the domain decomposition.

In the next section we investigate the numerical efficiency of the
proposed penalty parallel Robin--Robin type domain decomposition
schemes (\ref {eq_PII_85})~--~(\ref {eq_PII_86}) and (\ref
{eq_PII_91})~--~(\ref {eq_PII_92}).

\section{Numerical analysis}

We perform the numerical analysis of proposed domain decomposition
schemes for plane problems of unilateral contact between two
elastic bodies $\Omega_{\alpha} \subset {\mathbb R}^{2}$,
$\alpha=1,2$. For the numerical solution of linear variational
problems in subdomains, we use the finite element method (FEM)
with linear and quadratic triangular elements.

Note, that since proposed DDMs are obtained on the continuous
level, their convergence rate does not depend on the solution
methods in subdomains, if these methods are exact. If the methods,
used in subdomains are numerical, the number of iterations will
decrease as the precision of approximations in subdomains will
increase. Therefore, the proposed domain decomposition algorithms
are scalable.

At first, let us compare the convergence rates of different
particular domain decomposition schemes.

Consider the contact problem for two transversally isotropic
bodies $\Omega_{\alpha}$, $\alpha=1,2$, with the plane of
isotropy, parallel to the plane $x_{2}=0$ (Fig.~2) \cite {PII_8}.

The material properties of the bodies are: ${E_{\alpha}
\mathord{\left/ {\vphantom {E_{\alpha}  E'_{\alpha} }} \right.
\kern-\nulldelimiterspace} E'_{\alpha}} =2$, ${G_{\alpha}
\mathord{\left/ {\vphantom {G_{\alpha}  G'_{\alpha} }} \right.
\kern-\nulldelimiterspace} G'_{\alpha} } =2$, $\nu_{\alpha}
=\nu'_{\alpha} =0.3$, where $E_{\alpha}$, $\nu_{\alpha}$, and
$G_{\alpha}$ are the elasticity modulus, Poisson's ratio, and the
shear modulus for the body $\Omega_{\alpha}$ in the plane of
isotropy, and $E'_{\alpha}$, $\nu'_{\alpha}$, $G'_{\alpha}$ are
these constants in the orthogonal direction, $\alpha=1,2$.

The length and the height of each body is the same, and is equal
to $4\, b$. The distance between the bodies before the deformation
is $d_{12} ({\bf x})=r\, {x_{1}^{2} \mathord{\left/ {\vphantom
{x_{1}^{2}  b^{2} }} \right. \kern-\nulldelimiterspace} b^{2} }$,
the compression of the bodies is $\Delta \approx 2.154434\, r$,
$r=10^{-3} b$, and the possible contact area is $S_{12}
=\left\{{\bf x}=\left(x_{1} ,x_{2} \right)^{{\rm T}} :\, \, \, \,
x_{1} \in [0;\, \, 2\, b],\, \, x_{2} =4\, b\right\}$.

\begin{figure}[h]
 \center
{
\includegraphics[bb=0mm 0mm 208mm 296mm, width=60.1mm, height=80.9mm,
 viewport=3mm 4mm 205mm 292mm]{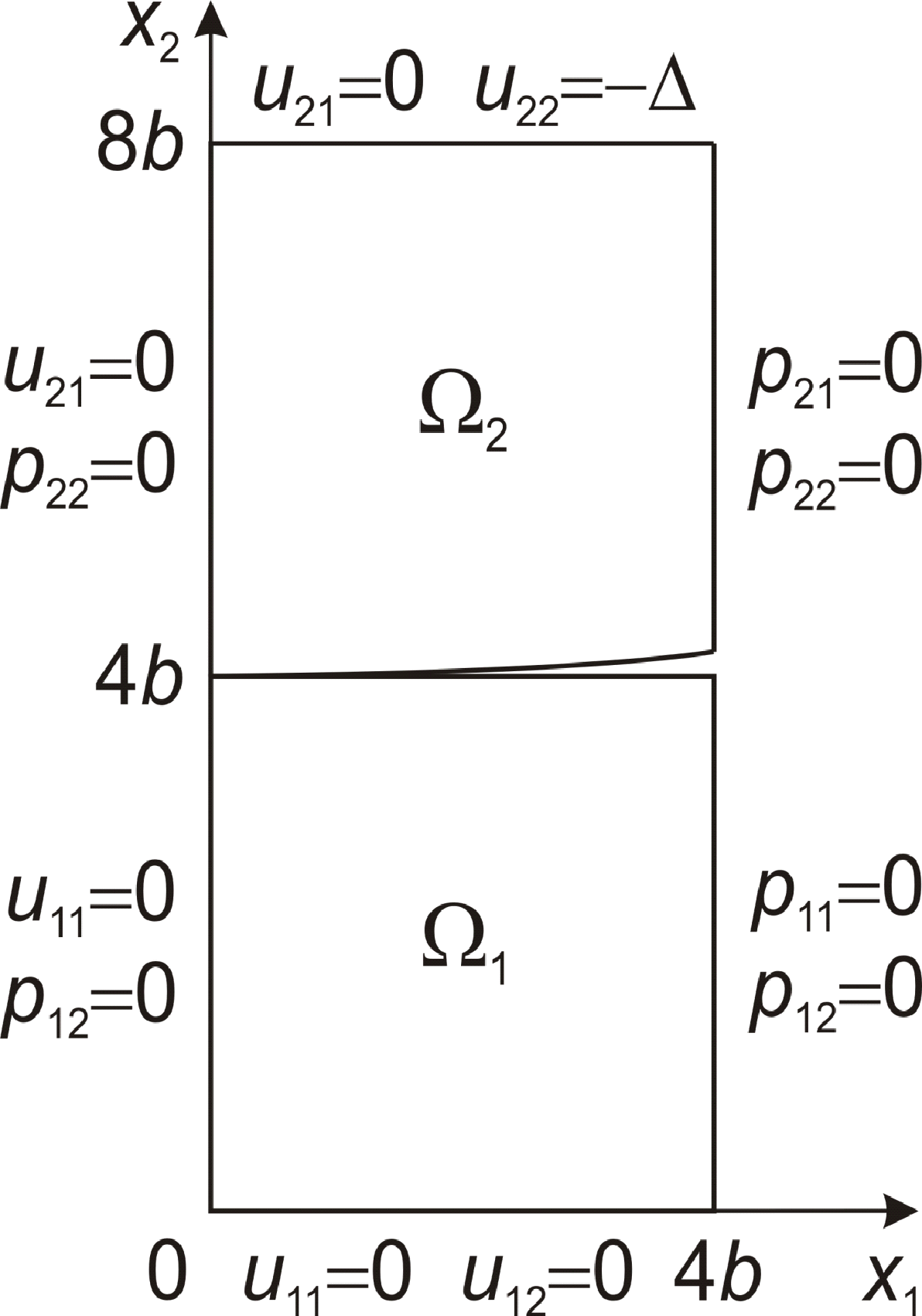}
 \\
 \textbf{Fig.~2.} Unilateral contact between two transversally isotropic bodies
}
\end{figure}
\begin{figure}[h]
 \center
{
\includegraphics[bb=0mm 0mm 208mm 296mm, width=82.3mm, height=79.4mm,
 viewport=3mm 4mm 205mm 292mm]{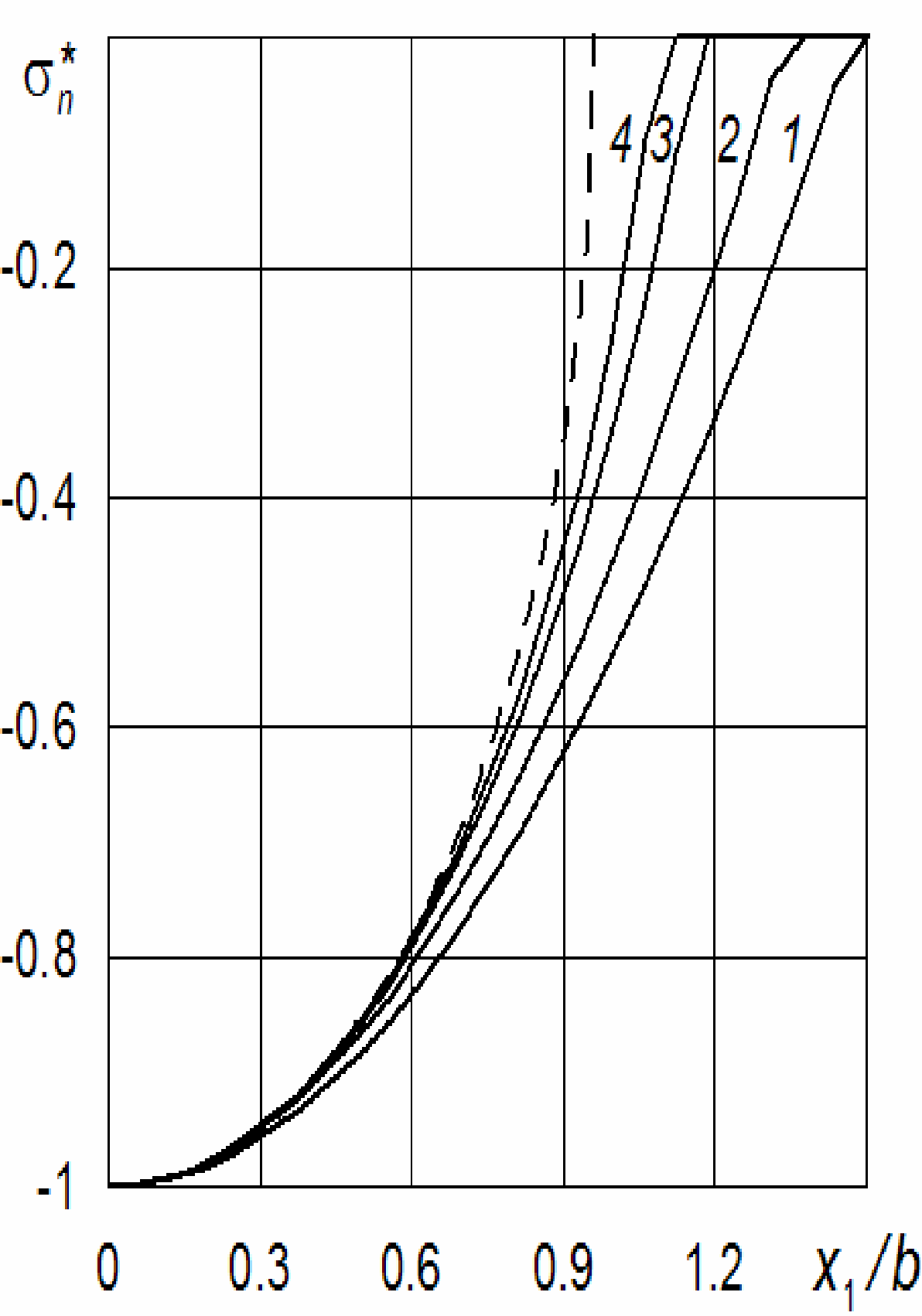}
\\
\textbf{Fig.~3.} Dimensionless normal contact stress at different
iterations
 }
\end{figure}

The problem was solved by parallel Robin--Robin domain
decomposition schemes, using FEM with 3190 quadratic triangular
elements in each body.

We took the penalty parameter in the form $\theta =4\,
bc\left({1\mathord{\left/ {\vphantom {1 E'_{1} }} \right.
\kern-\nulldelimiterspace} E'_{1} } +{1\mathord{\left/ {\vphantom
{1 E'_{2} }} \right. \kern-\nulldelimiterspace} E'_{2} } \right)$,
$c=0.05$, where $c$ is the dimensionless penalty coefficient, and
we used the following termination criterion for the domain
decomposition schemes:

\begin{equation} \label{eq_PII_103}
{\left\| u_{\alpha \, n}^{k+1} -u_{\alpha \, n}^{k} \right\| _{2}
\mathord{\left/ {\vphantom {\left\| u_{\alpha \, n}^{k+1}
-u_{\alpha \, n}^{k} \right\| _{2}  \left\| u_{\alpha \, n}^{k+1}
\right\| _{2} }} \right. \kern-\nulldelimiterspace} \left\|
u_{\alpha \, n}^{k+1} \right\| _{2} } \le \varepsilon _{u},\,\,\,
\alpha =1,2,...,N,
\end{equation}
where $\left\| u_{\alpha \, n} \right\| _{2} =\sqrt{\sum
_{j}\left[u_{\alpha \, n} ({\bf x}^{j} )\right]^{2}}$ is the
discrete norm, ${\bf x}^{j} \in S_{12}$ are the finite element
nodes on the possible contact area, and $\varepsilon_{u} > 0$ is
the relative accuracy for the displacements.

Fig.~3 shows the approximations of the dimensionless normal
contact stress $\sigma _{n}^{*} (x_{1} ,x_{2} )={\sigma _{12\,n}
(x_{1} ,x_{2} )\mathord{\left/ {\vphantom {\sigma _{12\,n} (x_{1}
,x_{2} ) \left|\sigma _{12n} (0,x_{2} )\right|}} \right.
\kern-\nulldelimiterspace} \left|\sigma _{12n} (0,x_{2}
)\right|}$, $x_{2} =4\, b$, $\left(x_{1} ,x_{2} \right)^{{\rm T}}
\in S_{12}$, obtained by the parallel Neumann--Neumann scheme
(\ref {eq_PII_88})~--~(\ref {eq_PII_89}) ($\psi _{12} ({\bf
x})=\psi _{21} ({\bf x})\equiv 0$) at iterations $k=1,\, \, 2,\,
\, 4,\, \, 21$ (Curves~1--4) for the optimal iterative parameter
$\bar{\gamma }=0.173$ and the accuracy $\varepsilon _{u}
=10^{-3}$. The dashed curve represents the exact solution for two
half-spaces, obtained in \cite {PII_27}. Hence, the real contact
area is $S_{12}^{*} \approx \left[\kern-0.15em\left[0;\, \,
b\right]\kern-0.15em\right]$, where $\left[\kern-0.15em\left[y;\,
\, z\right]\kern-0.15em\right]=\left\{{\bf x}=\left(x_{1} ,x_{2}
\right)^{{\rm T}} : \,\,\,\, x_{1} \in [y;\, \, z],\, \, x_{2}
=4\, b\right\}$.

At Fig.~4 and Fig.~5 the convergence rates of different particular
domain decomposition schemes are compared.

The dependence of the total number of iterations $m$ on the
iterative parameter $\gamma$ for the accuracy $\varepsilon_{u}
=10^{-3}$ is shown at Fig.~4, and its dependence on the
logarithmic accuracy $\lg \varepsilon_{u}$ for the optimal
iteration parameter $\gamma =\bar{\gamma}$ is shown at Fig.~5.

\begin{figure}[h]
 \center
{
\includegraphics[bb=0mm 0mm 208mm 296mm, width=79.1mm, height=76.6mm,
 viewport=3mm 4mm 205mm 292mm]{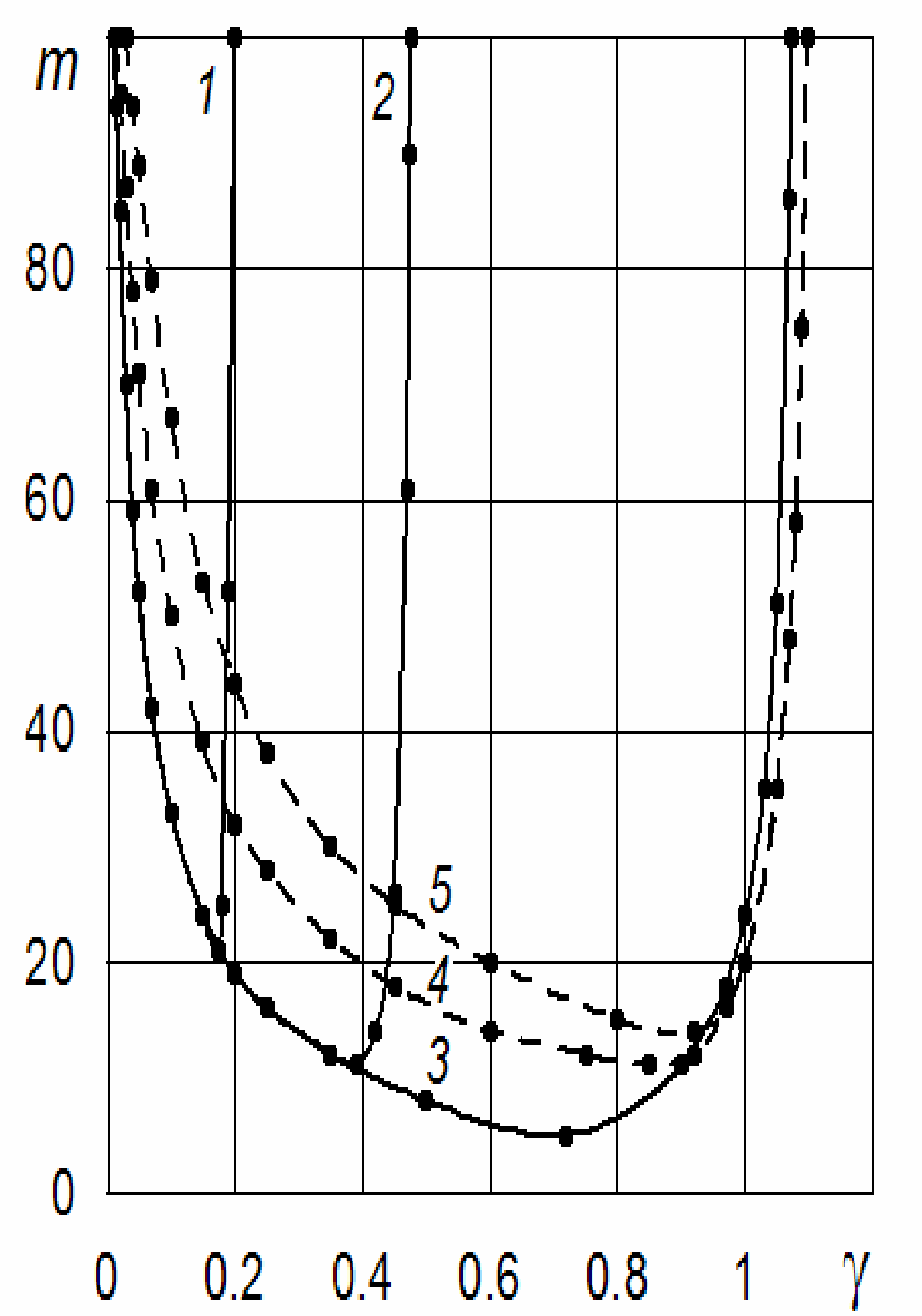}
 \\
 \textbf{Fig.~4.} The dependence of the total number of iterations on the iterative parameter
 $\gamma$
}
\end{figure}
\begin{figure}[h]
 \center
{
\includegraphics[bb=0mm 0mm 208mm 296mm, width=76.9mm, height=77.6mm,
viewport=3mm 4mm 205mm 292mm]{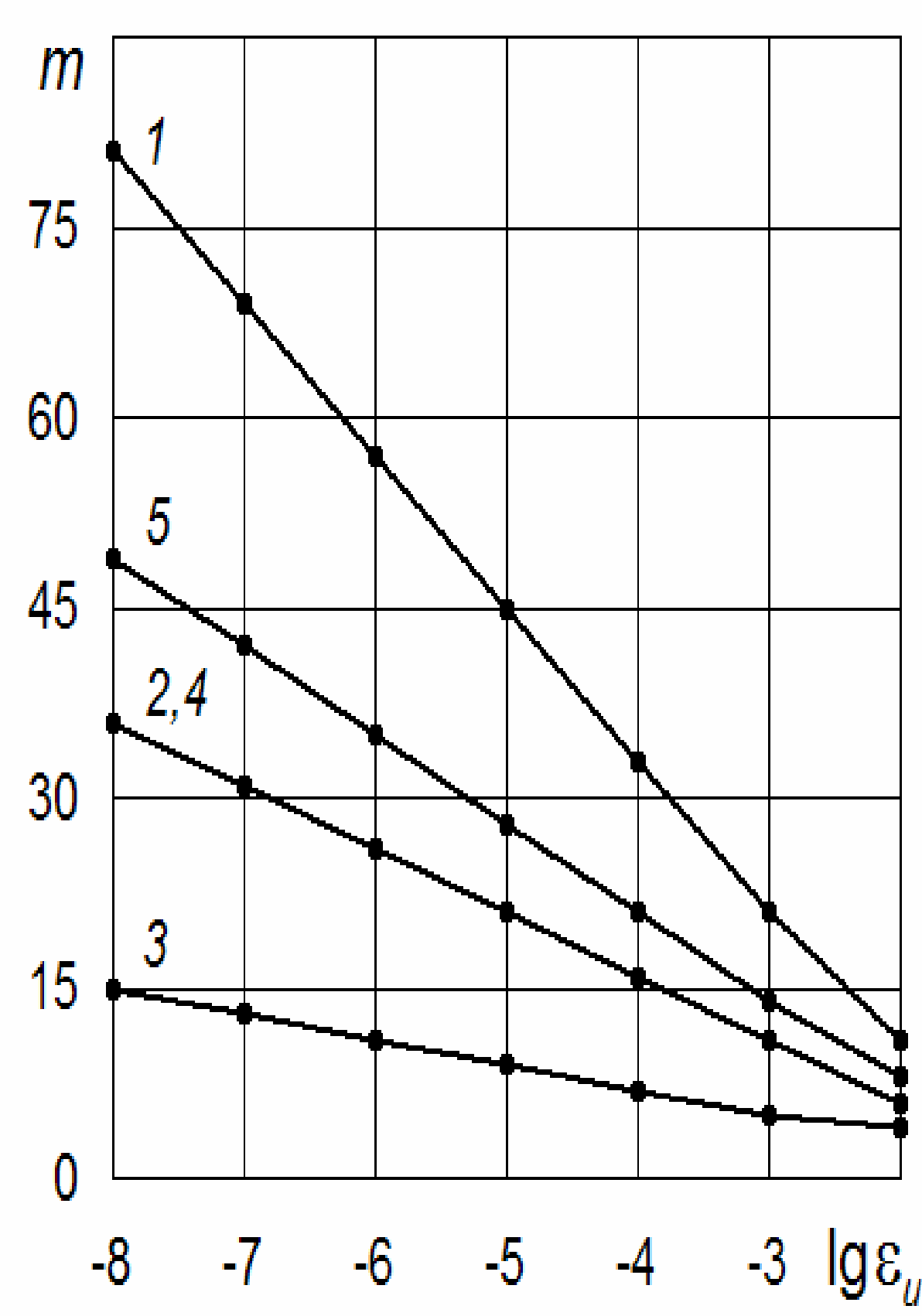}
\\
\textbf{Fig.~5.} The dependence of the total number of iterations
on the logarithmic accuracy
 }
\end{figure}

The first curve at these figures represents the parallel
Neumann--Neumann scheme ($S_{12}^{1} =S_{21}^{1} =\emptyset$,
$\psi _{12} ({\bf x})=\psi _{21} ({\bf x})\equiv 0$), Curves 2, 3,
4 and 5 correspond to the parallel Robin--Robin schemes (\ref
{eq_PII_85})~--~(\ref {eq_PII_86}) with $S_{12}^{1} =S_{21}^{1}$
equal to $\left[\kern-0.15em\left[0;\, \,
0.5\right]\kern-0.15em\right]$, $\left[\kern-0.15em\left[0;\, \,
1\right]\kern-0.15em\right]$, $\left[\kern-0.15em\left[0;\, \,
1.5\right]\kern-0.15em\right]$, and $\left[\kern-0.15em\left[0;\,
\, 2\right]\kern-0.15em\right]$ ($S_{12}^{1} =S_{21}^{1} =S_{12}
$, $\psi _{12} ({\bf x})=\psi _{21} ({\bf x})\equiv 1$)
respectively. Curve~3 also represents the nonstationary parallel
Dirichlet--Dirichlet scheme (\ref {eq_PII_97})~--~(\ref
{eq_PII_98}).

The optimal iterative parameter $\bar{\gamma}$ for the schemes
represented by Curves~1--5 is $\bar{\gamma}=0.173$, 0.39, 0.72,
0.85, and 0.92 respectively. For $\gamma =\bar{\gamma}$ and the
accuracy $\varepsilon_{u} =10^{-3}$ these schemes converge in 21,
11, 5, 11, and 14 iterations.

Thus, the convergence rate of the stationary Robin--Robin domain
decomposition schemes is linear. The parallel Robin--Robin scheme
(\ref {eq_PII_85})~--~(\ref {eq_PII_86}) with the surfaces
$S_{12}^{1}$, $S_{21}^{1}$ most closed to the real contact area
($S_{12}^{1} =S_{21}^{1} \approx S_{12}^{*} \approx
\left[\kern-0.15em\left[0;\, \, b\right]\kern-0.15em\right]$), and
the nonstationary parallel Dirichlet--Dirichlet scheme (\ref
{eq_PII_97})~--~(\ref {eq_PII_98}) ($\psi _{12} =\psi _{21} =\chi
_{12}^{k}$), which are represented by Curve~3, have the highest
convergence rates. These two schemes also have the widest range
from which the iterative parameter $\gamma $ can be chosen. The
convergence rate of the parallel Neumann--Neumann scheme
($S_{12}^{1} =S_{21}^{1} =\emptyset$), which is represented by
Curve~1, is the most slow.

Now let us investigate the convergence of the penalty method and
its dependence on the finite element discretization.

Consider the unilateral contact problem for two isotropic bodies
$\Omega_{1}$ and $\Omega_{2}$, one of which has a groove (Fig.~6).

The bodies are uniformly loaded by the normal stress with
intensity $q$. Each body has length $l$ and height $h$, and the
grove has length $b$.

\begin{figure}[h]
 \center
{
\includegraphics[bb=0mm 0mm 208mm 296mm, width=67.8mm, height=68.8mm,
 viewport=3mm 4mm 205mm 292mm]{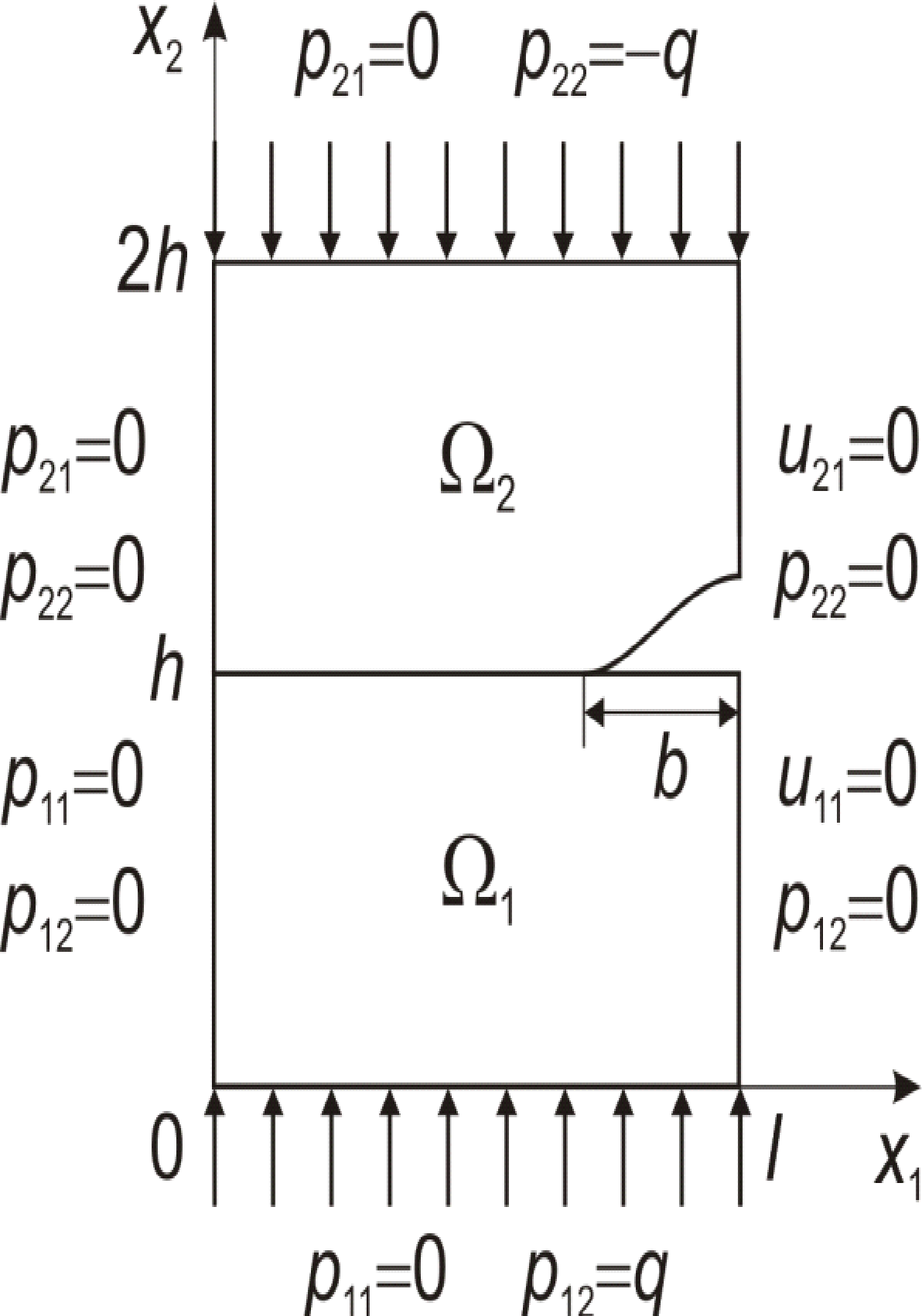}
 \\
 \textbf{Fig.~6.} Unilateral contact between two bodies with a groove
}
\end{figure}

\begin{figure}[h]
 \center
{
\includegraphics[bb=0mm 0mm 208mm 296mm, width=94.2mm, height=71.8mm,
 viewport=3mm 4mm 205mm 292mm]{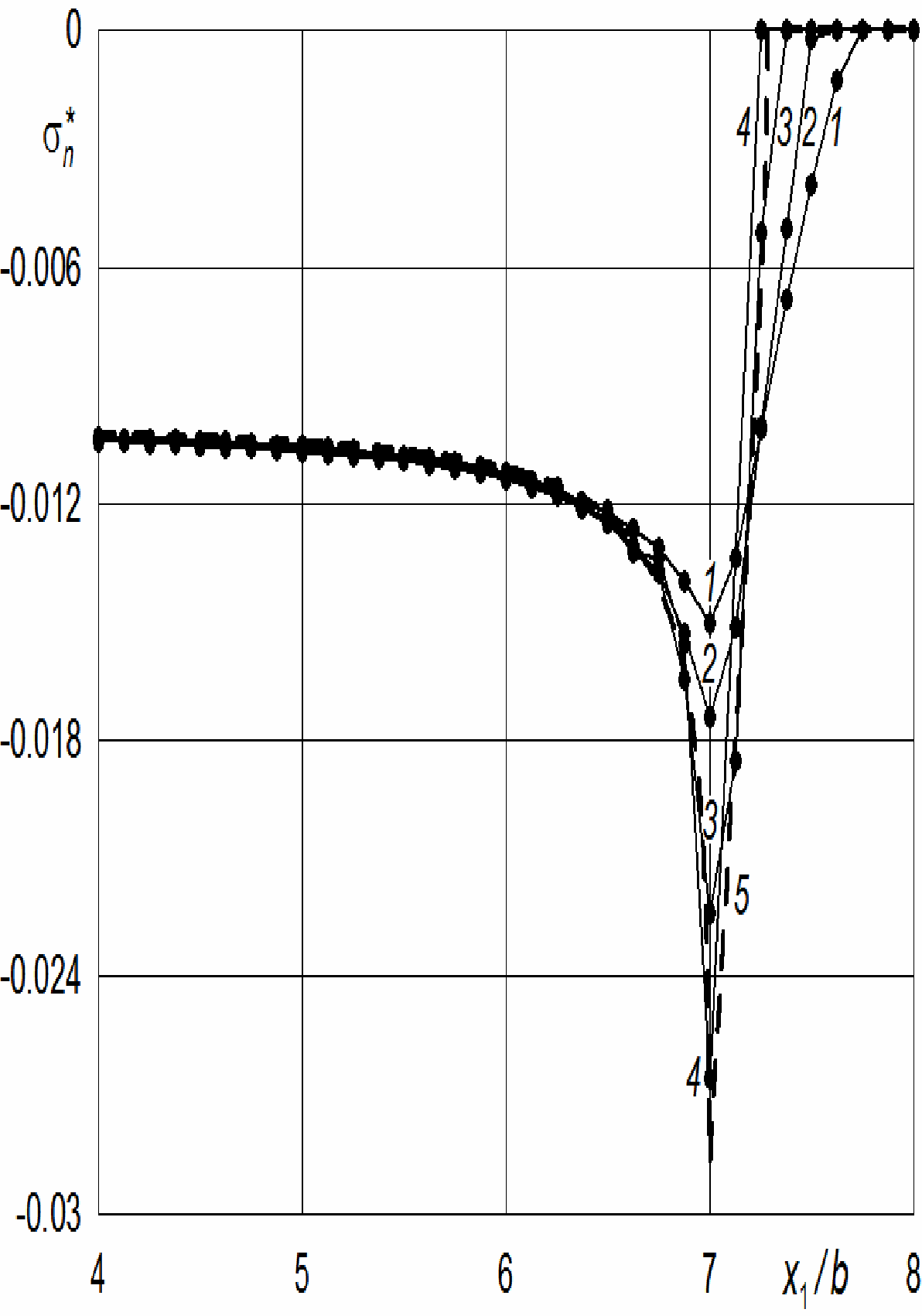}
 \\
\textbf{Fig.~7.} Normal contact stress $\sigma _{n}^{*} $ for
different penalty coefficients (at fixed finite element mesh) }
\end{figure}

The material properties of the bodies are the same: $E_{1} =E_{2}
=E$, $\nu _{1} =\nu _{2} =\nu =0.3$. The distance between the
bodies before the deformation is $d_{12} ({\bf x})=r\left\{{\, [\,
1-(x_{1} -l)^{2} \mathord{\left/ {\vphantom {\, [\, 1-(x_{1}
-l)^{2} b^{2} }} \right. \kern-\nulldelimiterspace} b^{2} } ]^{+}
\right\}^{3/2}$, where $r=0.05\, b$, $y^{+} =\max \{ 0,y\}$. The
possible contact area is $S_{12} =\left\{{\bf x}=\left(x_{1}
,x_{2} \right)^{{\rm T}} :\, \, \, \, x_{1} \in [0;\, \, l],\, \,
x_{2} =h\right\}$.

The exact solution of this problem in the case of the contact
between two half-spaces has a finite singularity in the flex point
of function $d_{12}({\bf x})$ \cite{PII_30}. Therefore, this
problem is a good test for the numerical methods.

The problem was solved by the nonstationary parallel
Dirichlet--Dirichlet domain decomposition scheme (\ref
{eq_PII_97})~--~(\ref {eq_PII_98}) with finite element
approximations on triangles.

The penalty parameter was taken as follows
\begin{equation} \label{eq_PII_104}
\theta =ch\, \sum _{\alpha =1}^{2}{\left(1-\nu _{\alpha }
\right)^{2} \mathord{\left/ {\vphantom {\left(1-\nu _{\alpha }
\right)^{2}  E_{\alpha } }} \right. \kern-\nulldelimiterspace}
E_{\alpha } },
\end{equation}
where $c$ is the dimensionless penalty coefficient. We used (\ref
{eq_PII_103}) as a termination criterion for the iterative
process.

For the iterative parameter $\gamma \in [0.45;\, \, 0.65]$ and the
accuracy $\varepsilon _{u} =10^{-3}$, and for the penalty
coefficients $c$ and the finite element meshes considered below,
the parallel Dirichlet--Dirichlet scheme, applied to solve this
problem, converges in 2--15 iterations.

Let us investigate the dependence of the quality of numerical
solution, obtained by this scheme, on the penalty parameter and
the finite element mesh.

Plots at Fig.~7 represent the approximations of the dimensionless
contact stress $\sigma _{n}^{*} ({\bf x})={\sigma _{12\,n} ({\bf
x})\mathord{\left/ {\vphantom {\sigma _{1\, 2n} ({\bf x}) E}}
\right. \kern-\nulldelimiterspace} E} $, ${\bf x} \in S_{12}$ for
the bodies with size $h=l=8\, b$ and external load $q=0.01\, E$,
obtained by the Dirichlet--Dirichlet scheme for different
dimensionless penalty coefficients $c$ at fixed finite element
mesh with 64 linear triangular finite elements on each side of the
possible contact area $S_{12}$. Curves 1--4 correspond to $c=0.1$,
$0.05$, $0.01$, and $0.0025$ respectively.

Plots at Fig.~8 represent the approximations of $\sigma _{n}^{*}
({\bf x})$, ${\bf x} \in S_{12}$ for the bodies with length $l=8\,
b$, height $h=2\, b$, and external load $q=0.0075\, E$, obtained
for different dimensionless penalty coefficients $c$ and different
finite element meshes. Curves 1 and 2 at this figure correspond to
$\sigma _{n}^{*}$ for the dimensionless penalty coefficients
$c=0.1$ and $c=0.01$ respectively at the finite element mesh with
32 linear triangular elements on each side of the possible
unilateral contact area $S_{12}$. Curves 3 and 4 correspond to
$\sigma _{n}^{*}$ for $c=0.1$ and $c=0.01$ respectively, but for
the finite element mesh with 64 linear triangular elements on each
side of $S_{12}$. Dashed curve at this figure and at Fig.~7
represents the exact solution, obtained in \cite {PII_30} for the
contact between two half-spaces.

\begin{figure}
 \center
{
\includegraphics[bb=0mm 0mm 208mm 296mm, width=95.1mm, height=69.1mm,
viewport=3mm 4mm 205mm 292mm]{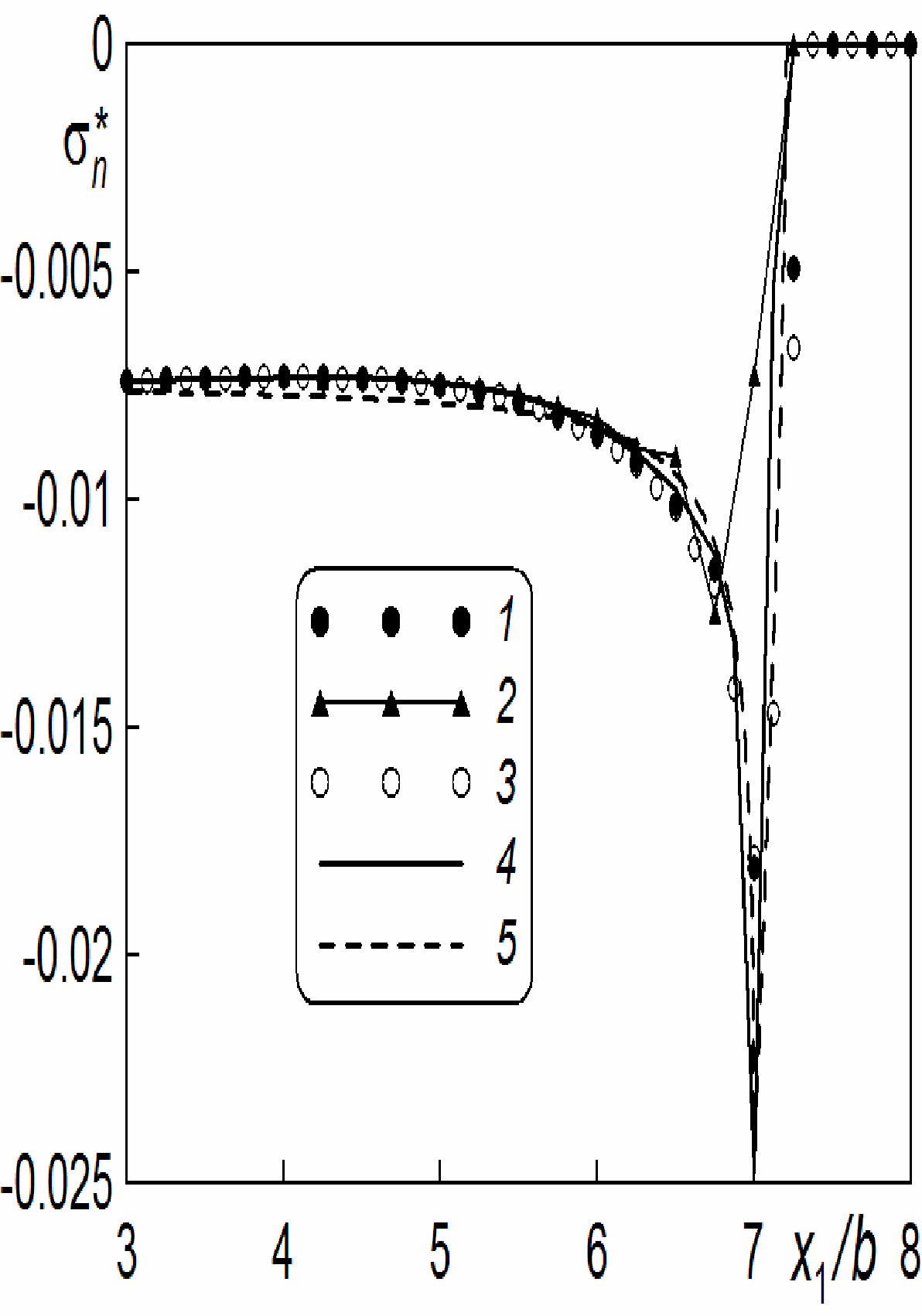}
 \\
 \textbf{Fig.~8.} Normal contact stress $\sigma _{n}^{*} $ for
different penalty coefficients and for different finite element
meshes}
\end{figure}

Here we see that in spite of the solution, obtained for the
penalty coefficient $c=0.1$ at the mesh with 32 finite elements on
each side of $S_{12}$ (Curve~1 at Fig.~8), the solution obtained
for the lower penalty coefficient $c=0.01$ at the same mesh
becomes instable. But if we refine the finite element mesh twice
for the penalty coefficient $c=0.01$, then the influence of the
errors on the perturbation of initial data will decrease, and we
will obtain much better approximation of the exact solution
(Curve~4 at Fig.~8).

Hence, we conclude that for obtaining a nice approximation of the
solution, we need to decrease the penalty parameter and to refine
the finite element mesh simultaneously.

\section{Conclusions}

For the solution of unilateral multibody contact problems of
elasticity we have proposed on the continuous level a class of
parallel Robin--Robin type domain decomposition schemes, which are
based on the penalty method for variational inequalities and some
stationary or nonstationary iterative methods for nonlinear
variational equations. In each iterative step of these schemes one
have to solve in parallel the linear variational equations in
subdomains, which correspond to some linear elasticity problems
with Robin boundary conditions on the possible contact areas.

We have given the mathematical justification of presented domain
decomposition methods. We have established the conditions of the
strong convergence of the solution of nonlinear penalty
variational equation, which corresponds to the original multibody
contact problem, to the weak solution of this problem.
Furthermore, we have proved theorems on the strong convergence and
stability of proposed DDMs, and have shown that the convergence
rate of stationary Robin--Robin schemes in some energy norm is
linear.

The numerical analysis of presented domain decomposition schemes
has been made for plane two-body contact problems using linear and
quadratic finite element approximations on triangles. The
convergence rates of different particular domain decomposition
schemes have been compared and their dependence on the iterative
parameter $\gamma$ has been investigated. The penalty parameter
and the mesh refinement influence on the numerical solution has
been examined. The numerical experiments have confirmed the
theoretical results on the convergence of these domain
decomposition schemes.

Among the positive features of proposed DDMs are the
regularization of the original contact problem because of the use
of the penalty term, the simplicity of their algorithms, and the
possibility to generalize them to more complicated contact
problems. Presented domain decomposition schemes allow to organize
parallel calculations and to use the most optimal mathematical
models (elastic body, shell theory) and discretization methods
(finite element method, boundary element method) for each of the
body (subdomain). These DDMs have only one iterative loop, which
deals simultaneously with the domain decomposition and the
nonlinearity of unilateral contact conditions. They do not require
to solve the nonlinear problems in each iterative step. Moreover,
since these methods are obtained on the continuous level, their
convergence rate does not depend highly on the discretization
techniques, i.e. the proposed domain decomposition algorithms are
scalable.


\newpage

\begin{thebibliography}{00}


\bibitem{PII_32}  P. Wriggers, Computational Contact Mechanics, second ed., Springer, Berlin Heidelberg, 2006.

\bibitem{PII_33}  B. Wohlmuth, Variationally consistent discretization schemes and numerical algorithms for contact problems, Acta Numerica \textbf{20} (2011) 569--734.

\bibitem{PII_31}  A. Toselli, O. B. Widlund, Domain Decomposition Methods -- Algorithms and Theory, Springer, Berlin Heidelberg, 2005.

\bibitem{PII_34}  A. Quarteroni, A. Valli, Domain Decomposition Methods for Partial Differential Equations, Oxford Science Publications, 1999.

\bibitem{PII_35}  B. F. Smith, P. E. Bj{\o}rstad, W. D. Gropp, Domain Decomposition, Cambridge Univ. Press, Cambridge, 1996.

\bibitem{PII_9}   I. Herrera, A. Carillo-Ledesma, A. Rosas-Medina, A brief overview of non-overlapping domain decomposition methods, Geof\'{i}sica Internacional \textbf{50}(4) (2011) 445--463.

\bibitem{PII_28}  Y. H. Savula, I. I. Dyyak, V. V. Krevs, Heterogeneous mathematical models in numerical analysis of structures, Computers and Mathematics with Applications \textbf{42}(8--9) (2001) 1201--1216.

\bibitem{PII_43}  P. L. Lions, On the Schwarz alternating method. III: A variant for nonoverlapping subdomains, Domain decomposition methods for partial differential equations, Proc. 3rd Int. Symp. Houston/TX (USA) 1989, 1990, pp. 202--223.

\bibitem{PII_44}  W. Guo, L. S. Hou, Generalizations and accelerations of Lions' nonoverlapping domain decomposition method for linear elliptic PDE, SIAM J. Numer. Anal. \textbf{41}(6) (2003) 2056--2080.

\bibitem{PII_45}  L. Qin, X. Xu, On  a  parallel  Robin-type  nonoverlapping domain decomposition  method, SIAM J. Numer. Anal. \textbf{44}(6) (2006) 2539--2558.

\bibitem{PII_36}  M. D. Gunzburger, J. S. Peterson, H. Kwon, An optimization based domain decomposition method for partial differential equations, Computers and Mathematics with Applications \textbf{37} (1999) 77--93.

\bibitem{PII_37}  M. D. Gunzburger, M. Heinkenschloss, H. K. Lee, Solution of elliptic partial differential equations by an optimization-based domain decomposition method, Applied Mathematics and Computation \textbf{113} (2000) 111--139.

\bibitem{PII_2}   G. Bayada, J. Sabil, T. Sassi, Algorithme de Neumann-Dirichlet pour des probl\`{e}mes de contact unilat\'{e}ral: R\'{e}sultat de convergence, C. R. Acad. Sci. Paris. Ser.~I~\textbf{335} (2002) 381--386 [In French].

\bibitem{PII_18}  R. Krause, B. A. Wohlmuth, Dirichlet-Neumann type algorithm for contact problems with friction, Computing and Visualization in Science \textbf{5}(3) (2002) 139--148.

\bibitem{PII_38}  C. Eck, B. Wohlmuth, Convergence of a Contact-Neumann iteration for the solution of two-body contact problems,  Mathematical Models and Methods in Applied Sciences \textbf{13}(8) (2003) 1103--1118.

\bibitem{PII_3}   G. Bayada, J. Sabil, T. Sassi, Neumann--Neumann domain decomposition algorithm for the Signorini problem, Appl. Math. Lett. \textbf{17}(10) (2004) 1153--1159.

\bibitem{PII_39}  J. Haslinger, R. Ku\v{c}era, T. Sassi, A domain decomposition algorithm for contact problems: Analysis and implementation, Math. Model. Nat. Phenom. \textbf{4}(1) (2009) 123--146.

\bibitem{PII_16}  J. Koko, An optimization-bazed domain decomposition method for a two-body contact problem, Num. Func. Anal. Optim. \textbf{24}(5--6) (2003) 586--605.

\bibitem{PII_40}  M. Ipopa, T. Sassi, Un algorithme de type Robin pour des probl\`{e}mes de contact unilat\'{e}ral, C. R. Acad. Sci. Paris. Ser.~I~\textbf{346} (2008) 357--362 [In French].

\bibitem{PII_26}  T. Sassi, M. Ipopa, F.-X. Roux, Generalization of Lions' nonoverlapping domain decomposition method for contact problems, Lect. Notes Comput. Sci. Eng. \textbf{60} (2008) 623--630.

\bibitem{PII_41}  M. Ipopa, T. Sassi, A Robin domain decomposition algorithm for contact problems: Convergence results, Lect. Notes Comput. Sci. Eng. \textbf{70} (2009) 145--152.

\bibitem{PII_17}  J. Koko, Uzawa block relaxation domain decomposition method for a two-body frictionless contact problem, Appl. Math. Lett. \textbf{22} (2009) 1534--1538.

\bibitem{PII_1}   P. Avery, C. Farhat, The FETI family of domain decomposition methods for inequality-constrained quadratic programming: Application to contact problems with conforming and nonconforming interfaces, Comput. Methods Appl. Mech. Engrg. \textbf{198} (2009) 1673--1683.

\bibitem{PII_5}   J. Dan\v{e}k, Domain decomposition method for contact problems with small range contact, J. Mathematics and Computers in Simulation \textbf{61}(3--6) (2003) 359--373.

\bibitem{PII_6}   Z. Dost\'{a}l, D. Hor\'{a}k, D. Stefanica, A scalable FETI--DP algorithm with non-penetration mortar conditions on contact interface, Journal of Computational and Applied Mathematics \textbf{231} (2009) 577--591.

\bibitem{PII_42}  Z. Dost\'{a}l, T. Kozubek, V. Vondr\'{a}k, T. Brzobohat\'{y}, A. Markopoulos, Scalable TFETI algorithm for the solution of multibody contact problems of elasticity, Int. J. Numer. Methods Eng. \textbf{41} (2010) 675--696.

\bibitem{PII_29}  J. Sch\"{o}berl, Efficient contact solvers based on domain decomposition techniques, Computers and Mathematics with Applications \textbf{42}(8--9) (2001) 1217--1228.

\bibitem{PII_23}  I. I. Prokopyshyn, Parallel domain decomposition schemes for frictionless contact problems of elasticity, Visnyk Lviv Univ., Ser. Appl. Math. Comp. Sci. \textbf{14} (2008) 123--133 [In Ukrainian].

\bibitem{PII_7}   I. I. Dyyak, I. I. Prokopyshyn, Convergence of the Neumann parallel scheme of the domain decomposition method for problems of frictionless contact between several elastic bodies, Journal of Mathematical Sciences \textbf{171}(4) (2010) 516--533.

\bibitem{PII_8}   I. I. Dyyak, I. I. Prokopyshyn, Domain decomposition schemes for frictionless multibody contact problems of elasticity, in: G. Kreiss et al (Eds.), Numerical Mathematics and Advanced Applications 2009. Proceedings of ENUMATH 2009, the 8th European Conference on Numerical Mathematics and Advanced Applications, Uppsala, July 2009, Springer, Berlin Heidelberg, 2010, pp. 297--305.

\bibitem{PII_47}  I. I. Prokopyshyn, I. I. Dyyak, R. M. Martynyak, I. A. Prokopyshyn,  Penalty Robin-Robin domain decomposition schemes for contact problems of nonlinear elasticity, Lect. Notes Comput. Sci. Eng. \textbf{91} (2013) 647--654.

\bibitem{PII_19}  A. S. Kravchuk, Formulation of the problem of contact between several deformable bodies as a nonlinear programming problem, Journal of Applied Mathematics and Mechanics \textbf{42}(3) (1978) 489--498.

\bibitem{PII_20}  V. I. Kuz'menko, On the variational method in the theory of contact problems for nonlinearly elastic laminated bodies, Journal of Applied Mathematics and Mechanics \textbf{43}(5) (1979) 961--970.

\bibitem{PII_4}  J. C\'{e}a, Optimisation: Th\'{e}orie et algorithmes, Dunod, Paris, 1971 [In French].

\bibitem{PII_21}  J.-L. Lions, Quelques m\'{e}thodes de r\'{e}solution des probl\`{e}mes aux limites non lin\'{e}aire, Dunod Gauthier-Villards, Paris, 1969 [In French].

\bibitem{PII_15}  N. Kikuchi, J. T. Oden, Contact Problem in Elasticity: A Study of Variational Inequalities and Finite Element Methods, SIAM, Philadelphia, 1988.

\bibitem{PII_22}  J.-L. Lions, E. Magenes, Probl\`{e}mes aux limites non homog\`{e}nes et applications, Volume~1, Dunod, Paris, 1968 [In French].

\bibitem{PII_11}  R. Glowinski, J.-L. Lions, R. Tr\'{e}moli\`{e}res, Analyse num\'{e}rique des in\'{e}quations variationnelles, Dunod, Paris, 1976 [In French].

\bibitem{PII_14}  A. M. Khludnev, V. A. Kovtunenko, Analysis of cracks in solids, WIT Press, Southampton, Boston, 2000.

\bibitem{PII_25}  I. I. Prokopyshyn, R. M. Martynyak, Numerical investigation of contact interaction of two solids with a groove by domain decomposition method, Problems of computational mechanics and strength of structures \textbf{16} (2011) 240--251 [In Ukrainian].

\bibitem{PII_12}  A. Ya. Grigorenko, I. I. Dyyak, S. I. Matysyak, I. I. Prokopyshyn, Domain decomposition methods applied to solve frictionless-contact problems for multilayer elastic bodies, Int. Appl. Mech. \textbf{46}(4) (2010) 388--399.

\bibitem{PII_48}  I. I. Prokopyshyn, I. I. Dyyak, R. M. Martynyak, I. A. Prokopyshyn,  Domain decomposition methods for problems of unilateral contact between elastic bodies with nonlinear Winkler covers, Lect. Notes Comput. Sci. Eng. \textbf{98} (2014) 739--748.

\bibitem{PII_49}  I. I. Prokopyshyn, Domain decomposition schemes based on penalty method for problems of ideal contact between elastic bodies, Mathematical methods and physicomechanical fields \textbf{57}(1) (2014) 41--56 [In Ukrainian].

\bibitem{PII_13}  A. Ya. Grigorenko, I. I. Dyyak, I. I. Prokopyshyn, Domain decomposition method with hybrid approximations applied to solve problems of elasticity, Int. Appl. Mech. \textbf{44}(11) (2008) 1213--1222.

\bibitem{PII_10}  M. Hinterm\"{u}ller, V. A. Kovtunenko, K. Kunisch, Generalized Newton methods for crack problems with non-penetration condition, Numer. Methods Partial Differential Eq. \textbf{21}(3) (2005) 586--610.

\bibitem{PII_27}  G. N. Savin, Pressure of an absolutely rigid stamp on an elastic anisotropic medium, Doklady Akademii Nauk SSSR \textbf{6} (1939) [In Russian].

\bibitem{PII_30}  R. M. Shvets, R. M. Martynyak, A. A. Kryshtafovych, Discontinuous contact of an anisotropic half-plane and a rigid base with disturbed surface, Int. J. Engng. Sci. \textbf{34}(2) (1996) 183--200.


\end{thebibliography}
\end{document}